\documentclass[12pt]{amsart}
\usepackage[shortlabels]{enumitem}
\usepackage[dvipsnames]{xcolor}
\usepackage{amssymb}
\usepackage[normalem]{ulem}

\usepackage{tikz}
\usepackage{pgfplots}
\usetikzlibrary{arrows,shapes,positioning}
\graphicspath{{graphics/}}

\tikzset{
    cross/.pic = {
    \draw[rotate = 30] (-#1,0) -- (#1,0);
    \draw[rotate = 60] (0,-#1) -- (0, #1);
    }
}

\usepackage{stmaryrd}
\usepackage{hyperref}
\usetikzlibrary {arrows.meta}

\newtheorem{thm}{Theorem}

\newtheorem{lem}[thm]{Lemma}

\newtheorem*{thm*}{Theorem}

\newtheorem{example}[thm]{Example}
\newcommand{\Z}{\mathbb Z}
\newcommand{\N}{\mathbb N}
\newcommand{\R}{\mathbb R}

\newcommand{\cA}{\mathcal A}
\newcommand{\cL}{\mathcal L}

\newcommand{\cP}{\mathcal P}
\newcommand{\cH}{\mathcal H}
\newcommand{\cM}{\mathcal M}
\newcommand{\cO}{\mathcal O}
\newcommand{\cU}{\mathcal U}
\newcommand{\cW}{\mathcal W}
\newcommand{\dist}{{\rm dist}}

\newcommand{\invmeas}{\mathcal M_T(X)}

\title[Asymptotic behavior of the pressure]{Asymptotic behavior of the pressure function for H\"older potentials}

\author{Tamara Kucherenko}\address{Department of Mathematics,
The City College of New York, New York, NY, 10031, USA}\email{tkucherenko@ccny.cuny.edu}

\author{Anthony Quas}\address{ Department of Mathematics and Statistics, University of Victoria, Victoria, BC
Canada}\email{aquas@uvic.ca}

\thanks{T.K. is supported by grants from the Simons Foundation \#430032. }
\thanks{A.Q. is supported by a grant from NSERC}

\begin{document}

\begin{abstract}
We study the behavior of the pressure function for H\"{o}lder continuous potentials on mixing subshifts of finite type. The classical theory of thermodynamic formalism shows that such pressure functions are convex, analytic and have slant asymptotes. 
We provide a sharp exponential lower bound on how fast the pressure function approaches its asymptotes. As a counterpart, we also show that there is no corresponding upper bound by exhibiting systems for which the convergence is arbitrarily slow. However, we prove that the exponential upper bound still holds for a generic H\"{o}lder potential. In addition, we determine that the pressure function 
satisfies a coarse uniform convexity property. Asymptotic bounds and quantitative convexity estimates are the first additional general properties of the pressure function obtained in the settings of Bowen and Ruelle since their groundbreaking work more than 40 years ago.

\end{abstract}

\keywords{thermodynamic formalism, topological pressure, variational principle, equilibrium states, H\"{o}lder potentials, subshifts of finite type}
\subjclass[2000]{}
\maketitle

\section{Introduction}\label{sec:1}
To quote Ruelle, ``the main object of the thermodynamic formalism is to study the differentiability and analyticity properties of the function $P$ [topological pressure], and the structure of the equilibrium states and Gibbs states" \cite[Introduction]{Ru1}. Arguably the most cited result in this context
is that for mixing subshifts of finite type the pressure is real analytic on the space of H\"{o}lder continuous potentials and that for each such potential there is only one Gibbs state which is also the only equilibrium state. These statements served as catalysts for the growth of the ergodic theory of smooth hyperbolic systems starting with Anosov maps. In his breakthrough work \cite{Bo} Bowen applied Ruelle's analytic tool of transfer operators to Anosov diffeomorphisms (in fact, more general Axiom A systems) using Markov partitions and symbolic coding introduced by Sinai. This allowed for the properties of the pressure and Gibbs states on shift spaces to be carried over to differentiable systems, resulting in a description of the behaviour of Lebesgue-almost every orbit.


Shortly after the introduction of the thermodynamic formalism, 
its relationship to dimension theory was discovered, where the concept of the topological pressure once again played a central role. 
A highly influential result due to Bowen \cite{Bo2} and Ruelle \cite{Ru2} is that the Hausdorff dimension of Julia sets for conformal maps can be computed as the root of the pressure function of a certain potential. It was used, in particular, to establish the analyticity of the Hausdorff dimension as a function of the parameter in the interior of the main cardioid of the Mandelbrot set. Since then dimensional estimates were obtained for numerous invariant sets and measures \cite{P,Ba,PU}, the vast majority of which use a version of Bowen's pressure formula.

As part of dimension theory, multifractal analysis is concerned with the complexity of level sets of asymptotically defined quantities such as Birkhoff averages, Lyapunov exponents, and local entropies. Usually, the geometry of a level set is sufficiently complicated to necessitate tools such as Hausdorff dimension or topological entropy in order to describe its size and complexity. In most cases, the main technical device to identify the various multifractal spectra is the pressure function, see e.g. \cite{BSS,BPS,C}. Through this approach the dimension of a level set is evaluated by the entropy of a suitable invariant measure which is produced as an equilibrium state for the appropriate potential. For instance, the pressure function of the geometric potential contains information about the spectrum of the maximum Lyapunov exponent for geodesic flows on compact manifolds \cite{BG}.


Despite the fact that the pressure function has been used in applications more and more over time, the understanding of the behavior of the function itself has not gone beyond the general statements of analyticity, convexity and existence of asymptotes -- properties already known to Bowen and Ruelle in the 1970's. Analyticity is the strongest possible regularity condition for real-valued functions. In the present work we examine the other two properties of the pressure in the classical settings of H\"older potentials and mixing subshifts of finite type.  We are able to characterize the rate of convergence of the pressure function to its asymptote as well as strengthen the convexity statement.




Throughout the paper we assume that $\phi:X\to\R$ is a H\"older continuous potential associated with a mixing subshift of finite type $(X,T)$. The \emph{topological pressure} of $\phi$ can be defined via the
Variational Principle by
\begin{equation*}
  P_{\rm top}(\phi)=\sup\left\{h_{T}(\mu)+\int\phi\,d\mu \right\}
\end{equation*}
where the supremum is taken over the set of all $T$-invariant probability measures on $X$ and $h_{T}(\mu)$
denotes the measure-theoretic entropy of the measure $\mu$. The measures which realize the above supremum
are called the \emph{equilibrium states} of $\phi$. The terminology comes from statistical physics: the quantity $E_\phi=-(h_T(\mu)+\int\phi\,d\mu)$ represents the free energy of the system in state $\mu$ and the equilibrium is given by the states which minimize the free energy. We refer the reader to the monographs \cite{Bo, Ru1, Walters} for a detailed exposition.

We study the \emph{pressure function} of $\phi$, $p_{\phi}(t)=P_{\rm top}(t\phi)$, where $t$ is a real-valued parameter. In statistical physics this function is regarded as a tool to observe an evolution of a system depending on a continuous external factor. One common interpretation of the parameter $t$ is the inverse temperature of the system. Then the behavior of $p_{\phi}(t)$ when $t\to\infty$ is of significant interest, since it reveals certain changes within the system when the temperature is lowered to zero. It has been observed that on the microscopic level materials tend to be highly ordered at a low temperature, which mathematically means that corresponding equilibrium states should be supported on configurations of low complexity \cite{VanEnter}. A system at absolute zero temperature exists in its ground state, hence the limit points of equilibrium states as temperature approaches zero are termed the ground states of the system.
A long standing conjecture in ergodic theory (finally resolved in the affirmative by Contreras in 2016 \cite{Co}) states that for a generic H\"{o}lder potential on a subshift of finite type the ground state is unique and supported on a periodic orbit. The question we address here is how fast the energy level of a system can approach the energy of its ground state when the temperature is lowered to zero. This leads to the task of characterizing the asymptotic behavior of the pressure function.





 For each $t$ the potential $t\phi$ has a unique equilibrium state $\mu_t$. The accumulation points of the family $(\mu_t)$ as $t\to\infty$  are the ground states of $\phi$. If $\mu_t$ converges (in the weak$^*$-topology) then the limit is called the zero-temperature measure.  The matter of existence of such a measure 
 received considerable attention in the literature. In 2001 Contreras, Lopes and Thieullen \cite{CLT} established the
existence of the zero-temperature limit for a generic set of H\"{o}lder potentials.
Two years later, Bremont \cite{Bremont} proved that any locally constant potential admits a zero-temperature measure, which piqued the interest in the validity of the same statement for H\"{o}lder potentials.  In 2010, Chazottes and Hochman \cite{ChH}  effectively ended the discussion by constructing an example of a
Lipschitz continuous potential on a full shift such that the zero-temperature limit does not exist.



The asymptotic behaviour of the pressure function as $t\to\infty$ is a classical result in ergodic optimization and thermodynamic formalism. Although it is often regarded as folklore, an early statement appears in a widely-cited but 
unpublished work of Conze and Guivarc'h \cite[Theorem 3.2]{ConzeGuivarc'h}. Related formulations were subsequently used by Jenkinson in his work on ergodic optimization \cite{Jenkinson}, while a proof in a more general setting of higher-dimensional actions could be found in \cite[Theorem 6.2]{ChazottesKucherenkoQuas}. Since we will use this result repeatedly, we record here a version adapted to shifts of finite type and continuous potentials.

\begin{thm*}[Folklore theorem on pressure asymptotics]
\label{thm:ConzeGuivarc'h}
Let $T$ be a shift map on a shift of finite type $X$.
Let $\phi$ be a continuous function on $X$. 
Define quantities $a$ and $b$ by
$$
b=\max_\mu\int\phi\,d\mu,
$$
where the maximum is taken over the (weak$^*$-compact) set $\mathcal M_{\text{inv}}$ of 
$T$-invariant measures; and
$$
a=\sup\left\{h(\mu)\colon \mu\in\mathcal M_{\text{inv}}, \int \phi\,d\mu=b\right\}.
$$
Then $p_\phi(t)\to a+bt$ as $t\to\infty$.
\end{thm*}

By weak$^*$-compactness and weak$^*$-upper semicontinuity of the entropy, there is 
at least one $T$-invariant measure $\mu_\infty$ such that $a=h(\mu_\infty)$ and $b=\int \phi\,d\mu_\infty$,
so that $p_\phi(t)\ge a+bt$ by the Variational Principle.
\color{black}

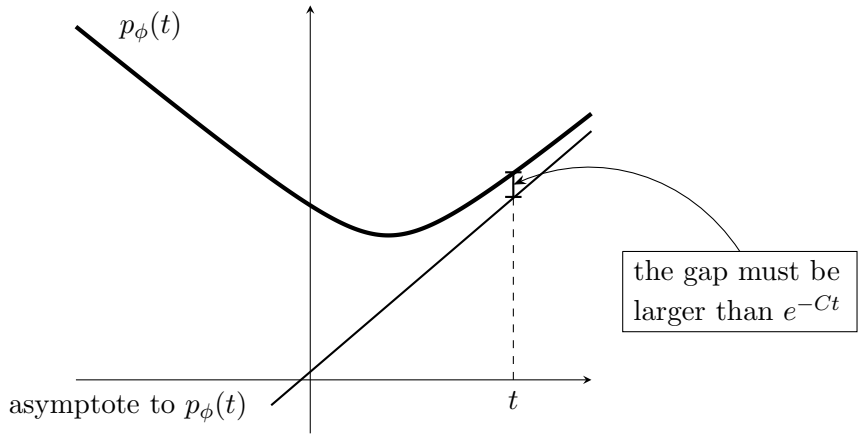
\begin{figure}[h]\label{pic:gap_at_infinity}
\[\begin{tikzpicture}
\begin{axis}[
samples=100,
axis y line=middle,
axis x line=middle,
ytick=\empty,
xtick=\empty,
clip=false,
xlabel near ticks,
ymax=7, ymin=-1,
]
\addplot[smooth,color=black,ultra thick,domain=-3:3.6]{sqrt((1.2*(x-1))^2+1)+1.7};
\addplot[smooth,color=black,thick, domain=-0.5:3.6]{1.25*(x-1)+1.4};
\node[left] at (axis cs:-1.5,6.6){\small{$p_\phi(t)$}};
\node[anchor=west] (source) at (axis cs:4.0,2){\small{the gap must be}};
\node[anchor=west] (source) at (axis cs:4.0,1.3){\small{$\text{larger than }e^{-Ct}$}};
\draw[-Stealth] (axis cs:5.5,2.4) to[bend left=-40] (axis cs:2.6,3.65);
\draw[dashed](axis cs:2.6,0) -- (axis cs:2.6,3.9);
    \node[anchor=north] at (axis cs:2.6,0){\small{$t$}};
\draw[|-|, thick](axis cs:2.6,3.4) -- (axis cs:2.6,3.9);
\draw(axis cs:4,2.4) -- (axis cs:7,2.4) -- (axis cs:7,0.9) -- (axis cs:4,0.9) -- (axis cs:4,2.4);
\node[anchor=west] (source) at (axis cs:-4.0,-0.5){\small{$\text{asymptote to }p_\phi(t)$}};
\end{axis}
\end{tikzpicture}\]
\caption{ This figure illustrates Theorem \ref{thm:gap_at_infinity}.}
\end{figure}

For any H\"older continuous $\phi$ that is not cohomologous to a constant, we establish a lower bound 
on the gap between the asymptote and the pressure function. We illustrate this statement in 
Figure \ref{pic:gap_at_infinity}.

\begin{thm}\label{thm:gap_at_infinity}
Let $X$ be a mixing subshift of finite type with positive entropy.
Let $\phi$ be a H\"older continuous function
that is not cohomologous to a constant.
Then there exist $C$ and $t_0$ such that $p_\phi(t)\ge\ell_\infty(t)+e^{-Ct}$
for all $t\ge t_0$, where $\ell_\infty(t)$ is the
asymptote to $p_\phi(t)$ at infinity.
\end{thm}

The exponential lower bound on the gap is the best one can hope for. We see in Example
\ref{ex:optimal_bnd} that for any locally constant potential on a full shift the rate with which
the pressure function approaches its asymptote is exactly exponential. This might suggest
that for a H\"{o}lder potential the upper bound on the rate should be of exponential type as well. However, this is
far off the mark. As it turns out, there is no upper bound at all as shown by the following theorem.

\begin{thm}\label{thm:large_gap}
Let $X$ be a mixing subshift of finite type with positive entropy and let $f:\R\to\R$ be any
convex function with an asymptote at infinity $\ell_\infty(t)=at+b$, where
$0\le b< h_{\rm top}(X)$.  Then there exists a H\"{o}lder potential $\phi:X\to \R$
such that $p_{\phi}(t)$ is asymptotic to $\ell_\infty(t)$ as $t\to\infty$ and $p_{\phi}(t)>f(t)$
for all sufficiently large $t$.
\end{thm}

Note that in order for the line $at+b$ to be an asymptote to the pressure function of some potential, $b$ must be the entropy of its ground state. Therefore, by the variational principle, $b$ cannot exceed the topological entropy of $X$ (or be negative). If $b=h_{\rm top}(X)$ then this ground state is necessarily the unique measure of maximal entropy of $X$, in which case the pressure function coincides with its asymptote. Hence, $0\le b< h_{\rm top}(X)$ is the weakest condition under which the statement holds.

We see from Theorem \ref{thm:large_gap} that the pressure function can converge to its asymptote arbitrarily slowly, while Theorem \ref{thm:gap_at_infinity} tells us that it cannot converge faster than exponentially. This raises natural questions: (i) what is the typical asymptotic behavior for H\"{o}lder potentials? (ii) do properties of the associated ground state have any impact on the rate of convergence? We provide a resolution to both. Question (ii) has a negative answer. In Section \ref{sec:genericity}  we construct two H\"{o}lder potentials on the full shift on two symbols which have the same zero temperature limit -- a point-mass measure. However, for one of them the pressure function approaches its asymptote exponentially fast, while for the other one the convergence is no faster than $\log\log t/\log t$. To answer Question (i), we prove that generically H\"{o}lder potentials admit an exponential upper bound on the gap between the pressure function and its asymptote. 

\begin{thm}\label{thm:gen_upper_bnd}
Let $X$ be a mixing shift of finite type and let $\cH$ denote a fixed H\"older class on $X$.
Then there is a dense open subset $\cU$ of potentials in $\cH$ such that for all $\phi\in \cU$,
there exist $C>0$ and $t_0$ such that $p_\phi(t)\le \ell_\infty(t)+e^{-Ct}$ for all $t\ge t_0$.
\end{thm}

This leads us to Question (iii): could a typical asymptotic behavior for H\"{o}lder potentials be given simply by $p_{\phi}(t)-\ell_\infty(t)\approx e^{-Ct}$ for some constant $C$?
We conjecture that the answer is yes, however the proof of such a statement appears to require an approach different from the one used in this work.

\begin{figure}[h]\label{pic:thm:strictconv}
\[\begin{tikzpicture}
\begin{axis}[
samples=100,
axis y line=middle,
axis x line=middle,
ytick=\empty,
xtick=\empty,
clip=false,
xlabel near ticks,
ymax=7, ymin=-1,
]
\addplot[smooth,color=black,ultra thick,domain=-3:3.6]{sqrt((1.2*(x-1))^2+1)+1.7};
\node[left] at (axis cs:-1.5,6.6){\small{$p_{\phi}(t)$}};
\node[anchor=west] (source) at (axis cs:4.0,2){\small{the gap must be}};
\node[anchor=west] (source) at (axis cs:4.0,1.3){\small{$\text{larger than }e^{-c(1+t^2)/h}$}};
\draw[-Stealth] (axis cs:5.85,2.4) to[bend left=-30] (axis cs:1.5,3.05);

\draw[dashed](axis cs:1.5,0) -- (axis cs:1.5,2.9);
    \node[anchor=north] at (axis cs:1.5,0){\small{$t$}};
\draw[dashed](axis cs:0.6,0) -- (axis cs:0.6,2.85);
    \node[anchor=north] at (axis cs:0.6,0){\small{$t-h$}};
\draw[dashed](axis cs:2.4,0) -- (axis cs:2.4,3.7);
    \node[anchor=north] at (axis cs:2.45,0){\small{$t+h$}};
\draw[|-|, thick](axis cs:1.5,2.85) -- (axis cs:1.5,3.3);
\draw(axis cs:0.6,2.85) -- (axis cs:2.4,3.7);

\draw(axis cs:4,2.4) -- (axis cs:7.85,2.4) -- (axis cs:7.85,0.9) -- (axis cs:4,0.9) -- (axis cs:4,2.4);
\end{axis}
\end{tikzpicture}\]
\caption{ This figure illustrates Theorem \ref{thm:strictconv}.}
\end{figure}
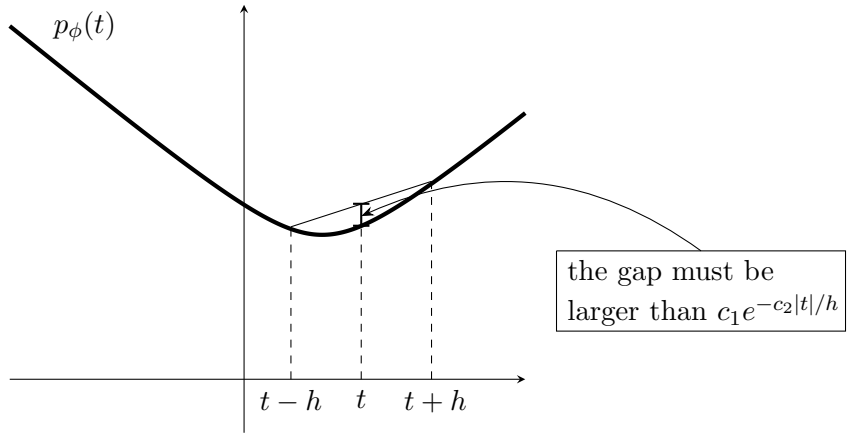

Next, we turn our attention to the convexity of the pressure function.
Although strictly convex analytic functions could be ``almost flat" on some intervals, we show that this is not possible for the pressure function of a H\"{o}lder potential. For a fixed $t\in\R$ we consider a symmetric interval $(t-h,t+h)$, where $h>0$. Since the pressure function is strictly convex, the midpoint of the secant line of the graph of the pressure function corresponding to points $t-h$ and $t+h$ is above the value of the pressure at $t$ (see Figure \ref{pic:thm:strictconv}). We show that the difference cannot be smaller than $e^{-c(1+t^2)/h}$ for some fixed positive constant $c$ which does not depend on the point $t$. We interpret this as a quantitative lower bound on the curvature, where the curvature bounds improve (a lot) when
one considers coarse intervals.

\begin{thm}\label{thm:strictconv}
Let $X$ be any mixing subshift of finite type with positive entropy and let $\phi$ be a H\"older continuous function that
is not cohomologous to a constant. Then there exists $c>0$
such that for any $t\in\R$ and any $h\in (0,1)$
\begin{equation}
\frac{p_{\phi}(t+h)+p_{\phi}(t-h)}{2}-p_{\phi}(t)>e^{-c(1+t^2)/h}.\label{eq:coarse}
\end{equation}
\end{thm}

One motivation for the present study comes from our previous work \cite{KQ} where we consider the pressure function of a \emph{continuous} potential on a full shift over a finite alphabet. It is known that such a pressure function is Lipschitz, convex, and has an asymptote at infinity. It turns out that in the case of continuous potentials no additional properties are present. In \cite{KQ} we explicitly construct a continuous potential on a full shift whose pressure function coincides with \emph{any} prescribed convex Lipschitz asymptotically linear function starting at a given positive value of the parameter. Immediately the question arose whether an analogously strong statement holds for the pressure function of a H\"{o}lder potential, where ``Lipschitz" is replaced by ``analytic". It follows from the assertions we made above that the answer is no. While this paper was in preparation we learned that the above question was also addressed in \cite{MaPollicott}, where a negative answer was obtained by establishing an inequality involving powers of the second, third, and fourth derivatives of the pressure function.


We finish this section with a brief outline of the paper. After preliminary material is given in Section \ref{sec:prelim}, the next two sections are devoted to the proofs of Theorem \ref{thm:gap_at_infinity} and Theorem \ref{thm:strictconv}.
 It is convenient to write $g_t(s)$ for the gap between $p_{\phi}(t)$ and its tangent line at $t$: $g_t(s)=
p_\phi(s)-\Big(h_T(\mu_t)+s\int\phi\,d\mu_t\Big)$ where $\mu_t$ is the (unique) equilibrium state for the potential
$t\phi$. We extend this notation to $g_\infty(s)$ for the gap between $p_\phi$ and its slant asymptote $\ell_\infty(s)=h_{T}(\mu_\infty)+s\int\phi\,d\mu_\infty$.
Recall that even though $\mu_\infty$ may not be unique, $\int\phi\,d\mu_{\infty}$ and $h_T(\mu_\infty)$ are the same for all
accumulation points, $\mu_\infty$, of the corresponding family of equilibrium states.
In this notation, Theorem \ref{thm:gap_at_infinity} may be re-expressed as
$$
g_\infty(t)\ge e^{-Ct}
$$
for all large $t$; and Theorem \ref{thm:strictconv} may be re-expressed as
$$
g_{t}(t-h)+g_{t}(t+h)\ge e^{-c(1+t^2)/h}
$$
for all $t\in\R$ and all $h>0$.

To prove the two theorems we estimate the gap functions $g_\infty$ and $g_t$ from below by building a new invariant measure $\mu'$ starting from $\mu_\infty$ (for Theorem \ref{thm:gap_at_infinity}) or $\mu_{t}$ (for Theorem \ref{thm:strictconv}) using
``coupling and splicing" techniques described in \cite{AQcutsplice}. A general construction of this kind is carried out in Section \ref{sec:CoupleSplice}. The objective is to increase the entropy of $\mu'$ compared to $\mu_\infty$(or $\mu_t$) while controlling the decrease in the value of the integral of $\phi$. In Section \ref{sec:main_proofs} we verify that the measure $\mu'$ can be constructed in such a way that the gain in entropy exceeds the drop in the integral ensuring that the quantity $h_T(\mu')+t\int\phi\,d\mu'$ is above the asymptote (or the tangent line) to the pressure function by the required amount.

In Section \ref{sec:slow_conv} we prove Theorem \ref{thm:large_gap}. We show that one can find a potential $\phi$
whose pressure function approaches its asymptote as slowly as desired. The idea of the proof is to define 
$\phi$ in terms of the distance to a carefully chosen subshift of $X$. We mimic the Rothstein 
shift \cite{Rothstein}, which has zero topological entropy,
but whose complexity can be chosen to exceed any given sub-exponential function.

Lastly, in Section \ref{sec:genericity} we establish the generic exponential upper bound for 
H\"{o}lder potentials (Theorem \ref{thm:gen_upper_bnd}). To prove this result we use the fact due to 
Contreras \cite{Co} that the set of potentials for which the zero temperature measure supported on a 
periodic orbit contains an open and dense set. However, this set must be necessarily modified for our 
purpose since there are potentials for which the zero temperature measure is supported on a periodic 
orbit, but for which the convergence of the pressure function to the asymptote is much slower than 
exponential, see Example \ref{ex:no_upper_bnd}.

The exposition in this article was greatly improved by following the suggestions
of the referees. We thank them for very helpful comments on the first version of this article. 

\section{Notation and conventions}\label{sec:prelim}
Our analysis takes place in the setting of two-sided shifts of finite type, which we always
assume to have forbidden blocks of length 2 only. We denote the alphabet by $\cA$
and the shift map by $T$. We use the metric $d(x,y)=2^{-n(x,y)}$
where $n(x,y)=\inf\{|n|\colon x_n\ne y_n\}$.


For a word $w=w_0\cdots w_{n-1}\in \cA^n$,
the cylinder generated by $w$ is denoted $[w]=\{x\in \cA^\Z: x_i=w_i\text{ for }$i$=0,\ldots,n-1\}$.
We write $|w|$ for its length, $n$ (and also refer to $n$ as the length of the cylinder set $[w]$).
Given any two words $w$ and $v$ we write
$wv$ for their concatenation. 
For a subshift $X\subset \cA^\Z$ we denote by $\cL_n(X)$
the set of all admissible words in $X$ of length $n$ and write
$\cL(X)=\bigcup_{n=1}^\infty \cL_n(X)$ for the language of $X$.
A subshift of finite type $X$ is \emph{mixing} if there exists an integer $L$ such that for any two words $u,v\in\cL(X)$ and any $n\ge L$ there is a word $w\in\cL_n(X)$ such that $uwv\in\cL(X)$. In this case
we call the smallest such $L$ the \emph{mixing length} of $X$.

We write $\cP$ for the (generating) partition consisting of all cylinder sets of length 1 and write
$\cP_n$ for the partition consisting of all cylinder sets of length $n$. If $\mu$ is a $T$-invariant measure
we denote its entropy by $h_T(\mu):=\lim_{n\to\infty}\frac 1n H_\mu(\cP_n)$,
where as usual $H_\mu(\mathcal Q)$ denotes the entropy of the countable partition $\mathcal Q$ with respect to measure $\mu$
(and we use natural logarithms in the definition).

A function $\phi$ is H\"older continuous if there exist $c\ge 0$ and $0<\alpha<1$ such that
$|\phi(x)-\phi(y)|\le c\alpha^{n(x,y)}$ where $n(x,y)$ is as above.
For a shift of finite type $X$ and a H\"older continuous potential $\phi$,
our principal object of study is the function $p_\phi(\cdot)$ given by $p_\phi(t)=P_{\rm top}(t\phi)$, where
$P_{\rm top}(\psi)$ denotes the topological pressure, equal by the Variational Principle to $\sup_\mu\left(
h_T(\mu)+\int \psi\,d\mu\right)$, where the supremum is taken over the collection of $T$-invariant
probability measures on $X$, $\cM(X)$.

We now recall basic properties of the function $p_\phi$ which are easily deduced from the Variational Principle. One may check that $p_\phi$ is convex. Monotonicity of the map $\psi\mapsto P_{\rm top}(\psi)$ together with the equality $P_{\rm top}(\psi+c)=P_{\rm top}(\psi)+c$ implies that $p_\phi$ is a Lipschitz function. 
As described in the folklore theorem on pressure asymptotics (Section \ref{sec:1}), $p_\phi$ has a slant asymptote given by $\ell_\infty(s)=h_T({\mu_\infty})
+s\int \phi\,d\mu_\infty$ where $\mu_\infty$ is any measure with maximal entropy among the
collection of maximizing measures.


Another property comes from the description of
equilibrium states as tangent functionals to the pressure given by Walters \cite{W1}.
If $\mu$ is any equilibrium state for $t\phi$, then the affine function
$\ell_\mu(s)=h_\mu(T)+s\int \phi\,d\mu$ is a sub-gradient of $p_\phi$ at $t$: $p_\phi(s)\ge \ell_\mu(s)$
for all $s$; and $\ell_\mu(t)=p_\phi(t)$. Conversely, for any sub-gradient $\ell(s)$ of
$p_\phi$ at $t$, there is an
equilibrium state $\mu$ such that $\ell(s)=\ell_\mu(s)$ for all $s$. Since $h_T(\mu)$ is a non-negative
quantity, bounded above by $\log |\cA|$,
we see that all subgradients of $p_\phi$ intercept the  vertical axis in a bounded sub-interval of $[0,\infty)$. One can show that boundedness of the vertical axis intercepts of the supporting lines implies both the
Lipschitz condition and the existence of a slant asymptote.


Although there are no other general properties of the pressure function for continuous potentials, one can say much more about the pressure under the restriction that $X$ is a subshift of finite type and potential $\phi:X\to\R$ is H\"{o}lder.
It was shown by Ruelle in \cite{Ru3} (see also \cite{Bo}) that in this case $\phi$ has a unique equilibrium state $\mu_\phi$ which satisfies \emph{Gibbs property}, namely there is a constant $C_\phi>0$ such that for all $n\in\N$, $w\in\cL_n(X)$  and $x\in [w]$ we have
\begin{equation}\label{eq:Gibbs}
\frac{1}{C_\phi}\le\frac{\mu_\phi([w])}{\exp\left(S_n\phi(x)-nP_{\rm top}(\phi)\right)}\le C_\phi,
\end{equation}
where $S_n\phi(x)$ denotes the Birkhoff sum $\sum_{i=0}^{n-1}\phi(T^i x)$. Throughout the paper, we take $C_\phi$ to be the optimal constant in \eqref{eq:Gibbs}.
Furthermore, the function $p_\phi(t)$ is analytic and strictly convex \cite{Ru1}.
Clearly in this case $p_\phi$ has a unique sub-gradient (in fact a tangent line)
at each $t$, which we denote by $\ell_t(\cdot)$, so that $\ell_t(s)=h(\mu_t)+s\int\phi\,d\mu_t$
where $\mu_t$ is the unique equilibrium state for $t\phi$.



\section{Construction of invariant measures}\label{sec:CoupleSplice}
In this section we use what have been termed ``coupling and splicing" techniques \cite{AQcutsplice} to build
a family of invariant measures $\mu'$ on a mixing subshift of finite type $(X,T)$ by modifying realizations of
an initial measure $\mu$. A word $w$ in the language of $X$ is fixed and
the measures $\mu'$ are, roughly speaking, obtained by starting from a realization $x$ of
$\mu$ and randomly replacing
words of length $|w|$ in $x$ by copies of $w$ with some
frequency. To ensure that the new point belongs to $X$, blocks of the mixing length $L$ prior to and
following the inserted $w$'s have to be modified also.

We start with some auxiliary measures that will be used in the construction.
We build a family of
ergodic measures on $\{0,1\}^\Z$ with two parameters $\eta$ and $M$ such that gaps between 1's are
independent and geometrically distributed with parameter $\eta$, 
taking values in $\{n\colon n\ge M\}$. To distinguish from the subshift $X$ we denote the shift map on $\{0,1\}^\Z$ by $\sigma$.

\begin{lem}\label{lem:build-nu}
There exists a family of ergodic invariant measures $\nu_{\eta,M}$ on $\{0,1\}^\Z$
where $\eta$ runs over $(0,1)$ and $M$ runs over $\N$ with the following
properties:
\begin{itemize}
\item (Spacing of 1's) For $\nu_{\eta,M}$-a.e.\ $y\in\{0,1\}^\Z$, if $i<j$ and $y_i=y_j=1$, then $j\ge i+M$;
\item (Frequency of 1's) 
$$
\nu_{\eta,M}([1])=\left(M+\frac{1}{e^{\eta}-1}\right)^{-1};
$$
\item (Entropy) 
$
h_\sigma(\nu_{\eta,M})
\ge \nu_{\eta,M}([1])\log\frac1\eta$. 
\end{itemize}
\end{lem}

\begin{proof}Set $S=\{M,M+1,M+2,\ldots\}^\Z$ with the shift map $\sigma$. For any parameter $\eta>0$,
equip $S$ with the Bernoulli measure $\chi$ where the symbol $M+k$ occurs with
probability $(1-e^{-\eta})e^{-k\eta}$ for $k=0,1,2,\ldots$. By a standard construction
we take the suspension of $S$ by the height function $h(x)=x_0$ to obtain the
space $\bar S=\{(x,n)\colon x\in S; 0\le n<x_0\}$, equipped with the standard
suspension map
$$
\bar\sigma(x,n)=\begin{cases}
(x,n+1)&\text{if $n<x_0-1$;}\\
(\sigma(x),0)&\text{if $n=x_0-1$.}
\end{cases}
$$
The measure $\chi$ on $S$ lifts to the probability measure $\bar\chi$ on $\bar S$
defined by
\begin{equation}\label{eq:lift}
 \begin{aligned}
\bar\chi(A\times \{j\})&=\chi(A)/\Big(\sum_{n\ge M}n(1-e^{-\eta})e^{-\eta(n-M)}\Big)\\
&=\chi(A)/\Big(M+(1-e^{-\eta})\sum_{n=0}^\infty ne^{-\eta n}\Big)\\
&=\chi(A)/\Big(M+1/(e^{\eta}-1)\Big),
\end{aligned}
\end{equation}
for any set $A\subset \{x\in S\colon x_0>j\}$. The denominator in this expression is simply
a normalization factor, which is just the integral of the height function. Since $\chi$ is an ergodic
$\sigma$-invariant measure, it follows that $\bar\chi$ is ergodic and $\bar\sigma$-invariant.

We then build a factor map from $\bar S$ to $Y=\{0,1\}^\Z$, defined by
\begin{equation}\label{eq:Phidef}
\Phi(x,n)_j=\begin{cases}
1&\text{if $\bar\sigma^j(x,n)\in\{(y,0)\colon y\in S\}$;}\\
0&\text{otherwise.}
\end{cases}
\end{equation}
The push-forward of $\bar\chi$ under $\Phi$ will be denoted by $\nu_{\eta,M}$,
or simply by $\nu$. The statement on the spacing of 1's is now immediate from the construction of $Y$. Taking $A=S$ and $j=0$ in (\ref{eq:lift}) gives the statement on the frequency of 1's.

Since $\Phi$ is one-to-one, the entropies $h_{\bar S}(\bar\chi)$ and $h_\sigma(\nu)$ are equal.
A calculation shows that
\begin{align*}
h_\sigma(\chi)&=-\sum_{n=0}^\infty(1-e^{-\eta})e^{-\eta n}\log \big((1-e^{-\eta})e^{-\eta n}\big)\\
&=-\log(1-e^{-\eta})+\eta(1-e^{-\eta})\sum_{n=0}^\infty ne^{-\eta n}\\
&=-\log(1-e^{-\eta})+\frac\eta {e^{\eta}-1}.
\end{align*}
Hence by Abramov's formula,
\begin{align*}
h_\sigma(\nu)&=h_{\bar \sigma}(\bar \chi)=
\nu([1])\left(-\log(1-e^{-\eta})+\frac\eta{e^{\eta}-1}\right)
\end{align*}
We verify that $-\log(1-e^{-\eta})+\frac\eta{e^\eta-1}\ge -\log\eta$ for $\eta\in (0,1)$ 
to obtain the estimate $h_\sigma(\nu_{\eta,M})
\ge \nu_{\eta,M}([1])\log\frac1\eta$.
\end{proof}

Now let $X$ be a shift of finite type (defined by forbidden blocks of length 2)
and $\mu$ be an ergodic invariant measure on $X$.
Let $w=w_0\ldots w_{m-1}$
be a word in the language of $X$ with the property that there does not exist $j<\frac {2m}3$ 
such that $w_0\ldots w_{n-1-j}=w_j\ldots w_{n-1}$ so that no two occurrences
of $w$ may be separated by less than $2n/3$.
We then say that $w$ satisfies the \emph{no long overlaps condition}. We now use the measures $\nu_{\eta,M}$ constructed above
to build a modified measure $\mu'_{\eta,w}$ on $X$ inserting additional $w$'s as described at the beginning of the section.

Let $M=m+2L+1$, where $L$ is the mixing length of $X$.
Then let $\nu_{\eta,M}$  and $(Y,\sigma)$ be as previously constructed. We build a new measure $\mu'_{\eta,w}$ as follows.
For each $a\in \cA$ denote by $u'(a)$ and $u''(a)$ the lexicographically minimal words of length $L$ such that $au'(a)w_0$ and $w_{m-1}u''(a)a$ belong to $\cL(X)$. Then for $a,b\in\cA$ let $v(a,b)=u'(a)wu''(b)$.
Define a map $\Phi\colon X\times Y\to X$ where $\Phi(x,y)_j$ is given by
\begin{equation*}
\begin{cases}v(x_{k-1},x_{k+2L+m})_{j-k}&\text{if $y_k=1$ for some $k\in\{j-m-2L-1,\ldots,j\}$};\\
x_j&\text{otherwise.}
\end{cases}
\end{equation*}

Informally, $\Phi(x,y)$ is the point $x'$ obtained by simultaneously replacing, for
each $k$ such that $y_k=1$, the word $x_{k+L}\ldots x_{k+L+m-1}$
by $w$, and choosing $x'_k\ldots x'_{k+L-1}$ and $x'_{k+L+m}\ldots x'_{k+m+2L-1}$
to be the minimal words so that the resulting word belongs to $X$. We then obtain a measure $\mu'_{\eta,w}$ on $X$ defined by
\begin{equation}\label{eq:CoupleSplice}
\mu'_{\eta,w}=\Phi_*(\mu\times\nu_{\eta,M}).
\end{equation}
We refer to this as the \emph{coupling and splicing construction} of the measure $\mu'_{\eta,w}$ from $\mu$. 



\begin{lem}\label{lem:ent_estimate}
Let $(X,T)$ be a mixing subshift of finite type with mixing length $L$,
$\mu$ be an ergodic invariant measure on $X$, and
$w$ be a word of length at least $3L$ with no long overlaps.
For a given $\eta\in(0,1)$, define
$$
\delta:=\nu_{\eta,2L+|w|+1}([1])
=\frac{1}{(2L+|w|+1)+1/(1-e^{-\eta})},
$$
as computed in Lemma~\ref{lem:build-nu}. If $\eta\ge\mu([w])$ then
$\delta\ge \frac{\mu([w])}{5}$,
and the invariant measure $\mu'_{\eta,w}$ constructed in
\eqref{eq:CoupleSplice} satisfies
\begin{equation}\label{eq:hbound}
h_T(\mu'_{\eta,w})
\ge
h_T(\mu)
+h_{\sigma}(\nu_{\eta,2L+|w|+1})
-\delta H_\mu(\mathcal P_{2L+|w|})
-8\delta\log 2 .
\end{equation}
\end{lem}

\begin{proof}
We write $\bar X$ for the space $X\times Y\times X$ and
$\bar\Phi\colon X\times Y\to \bar X$ for the map
$\bar\Phi(x,y)=(x,y,\Phi(x,y))$, where $\Phi$ is as defined in \eqref{eq:Phidef}.
Let $\nu=\nu_{\eta,2L+|w|+1}$, $\mu'=\mu'_{\eta,w}$, and $\bar\mu$ be the
measure $\bar\Phi_*(\mu\times\nu)$ on $\bar X$. The product of the
three shift maps (one in each coordinate) is denoted by $\bar T$, i.e. $\bar T=T\times\sigma\times T$.
We introduce three partitions of $\bar X$:
\begin{align*}
\mathcal P_X&=\Big\{\{(x,y,z):x_0=a\}: a\in \cA\Big\};\\
\mathcal P_Y&=\Big\{\{(x,y,z):y_0=\epsilon\}\colon \epsilon\in\{0,1\}\Big\};\text{ and}\\
\mathcal P_Z&=\Big\{\{(x,y,z):z_0=a\}\colon a\in \cA\Big\}.
\end{align*}
We have the following equalities:
\begin{align*}
h_T(\mu')&=h_{\bar T}(\bar\mu,\mathcal P_Z)\\
h_{\bar T}(\bar\mu)&=h_{T\times\sigma}(\mu\times\nu)=h_T(\mu)+h_\sigma(\nu)\\
h_{\bar T}(\bar\mu)&=h_{\bar T}(\bar\mu,\mathcal P_Z)+h_{\bar T}(\bar \mu,\mathcal P_Y|\mathcal F_Z)
+h_{\bar T}(\bar \mu,\mathcal P_X|\mathcal F_Y\vee\mathcal F_Z),
\end{align*}
where $\mathcal F_Z=\bigvee_{j=-\infty}^\infty \bar T^{-j}\mathcal\mathcal P_Z$ with a similar definition
for $\mathcal F_Y$; and where the second equality follows since $\bar\Phi$ is an isomorphism from
$X\times Y$ equipped with the measure $\mu\times\nu$ to $\bar X$ equipped with the measure $\bar\mu$.

Combining the equalities gives
\begin{equation}\label{eq:condent}
h_T(\mu')=h_T(\mu)+h_\sigma(\nu)-h_{\bar T}(\bar \mu,\mathcal P_Y|\mathcal F_Z)-
h_{\bar T}(\bar \mu,\mathcal P_X|\mathcal F_Y\vee\mathcal F_Z)
\end{equation}
Hence to obtain a lower bound for $h_T(\mu')$, we require upper bounds for
$h_{\bar T}(\bar \mu,\mathcal P_Y|\mathcal F_Z)$ and
$h_{\bar T}(\bar \mu,\mathcal P_X|\mathcal F_Y\vee\mathcal F_Z)$.

Let $\delta=\nu([1])$, so that by Lemma \ref{lem:build-nu}, using $M=2L+|w|+1$,
\begin{equation}\label{eq:delta}
\delta=\frac1{(2L+|w|+1)+1/(1-e^{-\eta})}.
\end{equation}
Recall that $|w|\ge 3L\ge 2L+1$. The function $\eta\mapsto 1-e^{-\eta}$ is concave, so the graph over the interval $[0,1]$ lies above the secant line $\eta\mapsto (1-e^{-1})\eta$ so that $(1-e^{-1})\eta\le 1-e^{-\eta}$. Since $(1-e^{-1})^{-1}<2$, $$\frac{1}{1-e^{-\eta}}\le\frac{1}{(1-e^{-1})\eta}<\frac{2}{\eta},$$ so that $$\delta>\frac{1}{2|w|+2/\eta}.$$
By the no long overlap condition, we see the sets $\sigma^i([w])$ are disjoint for 
$i=0,\ldots,\lceil\frac {2|w|}3\rceil-1$
so that we obtain the very weak (but general) bound $\mu([w])\le \frac{3}{2|w|}$. Since by assumption $\mu([w])\le\eta\le 1$, we see that $2|w|+2/\eta\le3/\mu([w])+2/\mu([w])=5/\mu([w])$ and  $\delta\ge \mu([w])/5$
as claimed.

We next claim that
$$
h_{\bar T}(\bar \mu,\mathcal P_X|\mathcal F_Y\vee\mathcal F_Z)
\le \delta H_{\mu}(\mathcal P_{2L+|w|}).
$$
This is very intuitive in terms of information: given the $y$ and $z$ strings, $x$ matches $z$ outside
blocks of length $2L+|w|$ starting at each $k$ where $y_k=1$. Since $y$ is independent of $x$,
the average amount of information in a single reconstruction is
$H_{\mu}(\mathcal P_{2L+|w|})$.
One then expects that $\delta H_{\mu}(\mathcal P_{2L+|w|})$
is an upper bound for $h_{\bar T}(\bar \mu,\mathcal P_X|\mathcal F_Y\vee\mathcal F_Z)$
since if the blocks are sequentially reconstructed, knowledge of previous reconstructions
may give you some information about the current reconstruction.
For a formal proof, consider the induced map on $E=\{(x,y,z)\in\bar X\colon y_0=1\}$ and note that
$\bar\mu(E)=\delta$.
Let $r_E(x,y,z)=\min\{n>0\colon y_n=1\}$.
We introduce countable partitions of $E$: $\tilde{\mathcal P}_X$ and $\tilde{\mathcal P}_Z$
where the elements of $\tilde{\mathcal P}_X$ are of the form
$$
\left\{E\cap B_{i_0}\cap \bar T^{-1}B_{i_1}\cap \ldots\cap \bar T^{-(n-1)}B_{i_{n-1}}
\cap \{\bar x\colon r_E(\bar x)=n\}\right\},
$$
where $n$ runs over $\mathbb N$ and $B_{i_0},\ldots,B_{i_{n-1}}$ are elements of $\mathcal P_X$.
We define $\tilde{\mathcal P}_Z$ analogously. Thus $\tilde{\mathcal P}_X$ and $\tilde{\mathcal P}_Z$
partition $\bar X$ according to the return time to $E$ and the symbols in the $x$- and $z$-coordinates
until that return respectively.

By Abramov's formula, we have
\begin{equation}\label{eq:Abramov}
h_{\bar T}(\bar \mu,\mathcal P_X|\mathcal F_Y\vee\mathcal F_Z)=
\delta h_{\bar T_E}(\bar\mu_E,\tilde{\mathcal P}_X|\tilde{\mathcal F}_Y\vee\tilde{\mathcal F}_Z).
\end{equation}
Since $E\in\mathcal F_Y\vee \mathcal F_Z$, $r_E$ is ($\mathcal F_Y\vee \mathcal F_Z$)-measurable.

As a partition of $E$, $(\mathcal P_X)_{2L+|w|}\vee \tilde{\mathcal P}_Z=
\tilde{\mathcal P}_X\vee\tilde{\mathcal P}_Z$ since $z_k$ agrees with $x_k$ everywhere except on
the $2L+|w|$ symbols following an occurrence of 1 in the $y$ coordinate.
That is $\mathcal F_Y\vee\mathcal F_Z\vee(\mathcal P_X)_{2L+|w|}$
 is a refinement of $\tilde{\mathcal P}_X$.
It follows that
\begin{align*}
h_{\bar T_E}(\bar\mu_E,\tilde{\mathcal P}_X|\tilde{\mathcal F}_Y\vee\tilde{\mathcal F}_Z)
&\le h_{\bar T_E}(\bar\mu_E,(\mathcal P_X)_{2L+|w|})\\
&\le
H_{\bar\mu_E}((\mathcal P_X)_{2L+|w|})\\
&=H_\mu(\mathcal P_{2L+|w|}),
\end{align*}
where for the last equality we used the independence of $\mathcal F_X$ and
$\mathcal F_Y$.
Combining this with \eqref{eq:Abramov}, we see
\begin{equation}\label{eq:h(X|YZ)}
h_{\bar T}(\bar \mu,\mathcal P_X|\mathcal F_Y\vee\mathcal F_Z)
\le \delta H_{\mu}(\mathcal P_{2L+|w|}).
\end{equation}


We next estimate $h_{\bar T}(\bar\mu,\mathcal P_Y|\mathcal F_Z)$. We use the hypothesis that
$\delta\ge \frac 15\mu([w])$. Given this, we need an estimate
for $\mu'([w])$. There are three ways that the word $w$ may appear in $z=\Phi(x,y)$: firstly if $y_{k-L}=1$, then
$z_k...z_{k+|w|-1}=w$; secondly if $x_k...x_{k+|w|-1}=w$ then $z_k...z_{k+|w|-1}$ may also be $w$ (unless
the map $\Phi$ overwrites some of that part of $x$); and thirdly if a $w$ is ``inadvertently" created involving
some parts of the original sequence $x$ and some symbols that are modified by the map $\Phi$. Since the instances of $w$ are at least $\frac 23 |w|$ apart by the no long overlaps condition, we claim there can be at most two of this third type of $w$ for each instance of the first type.
To see this, notice that the next possible $w$ occurs at least $\frac{2m}3$ positions to the right of $k+L$; and the next $w$ after that is at least $\frac{4m}3$ positions to the right of $k+L$. However since $L<\frac m3$, this second $w$ lies outside the coordinate range $k,\ldots, k+2L+m-1$ that is modified as a result of $y_k$ being equal to 1. This shows that for each $w$ inserted as a result of $y_k$ being equal to 1, there is at most one ``inadvertent" $w$ formed to the right of the inserted $w$. A similar argument shows that there is at most one inadvertent $w$ to the left of the inserted $w$.
We see that $\mu'([w])\le \mu([w])+3\delta\le 8\delta$.

Given this, we estimate $h_{\bar T}(\bar\mu,\mathcal P_Y|\mathcal F_Z)$
using the induced system of $\bar T$ with returns to $G:=\{(x,y,z):z_0...z_{|w|-1}=w\}$. Note that
$\bar\mu(G)=\mu'([w])$.
We then let $\mathcal Q_Y$ be the partition of $G$ according to the first return time map $r_G(\bar x)=\min\{n>0:\bar T^n(\bar x)\in G\}$, and certain symbols in the $y$ sequence.
Specifically $\mathcal Q_Y$ is the countable partition of $G$ with elements of the form
$$
G\cap \sigma^{-L}B_{i_0}\cap \sigma^{-L+1}B_{i_1}\cap\ldots\cap \sigma^{n-1-L}B_{i_{n-1}}\cap
\{\bar x\in G\colon r_G(\bar x)=n\},
$$
where $n$ runs over the positive integers and the $B$'s are elements of $\mathcal P_Y$. That is $G$ is partitioned according
to the return time and the $y$-symbols between time $-L$ and $r_G-L-1$. By the construction of $\bar \Phi$, if
$\bar x=(x,y,z)\in G$, $y_{-L}$ is either 1 or 0 (according to whether the copy of $w$ was deliberately inserted or not)
and $y_{-L+1},\ldots,y_{r_G(\bar x)-L-1}$ are all 0.
By Abramov's theorem,
$$
h_{\bar T}(\bar\mu,\mathcal P_Y|\mathcal F_Z)=\bar\mu(G)h_{\bar T_G}(\bar\mu_G,\mathcal Q_Y|\mathcal F_Z).
$$
By the above description, since $G$ is $\mathcal F_Z$-measurable,
$H_{\bar\mu_G}(\mathcal Q_Y|\mathcal F_Z)\le\log 2$, so that
\begin{equation}\label{eq:h(Y|Z)}
h_{\bar T}(\bar \mu,\mathcal P_Y|\mathcal F_Z)\le 8\delta\log 2.
\end{equation}

Substituting
\eqref{eq:h(X|YZ)} and \eqref{eq:h(Y|Z)} in \eqref{eq:condent} gives \eqref{eq:hbound} as required.
\end{proof}

\begin{lem}\label{lem:int_estimate}
Let $X$ be a mixing subshift of finite type with mixing length $L$ and $\mu$ be an 
ergodic measure on $X$. Suppose $\phi:X\to\R$ is H\"older, satisfying 
$|\phi(x)-\phi(y)|\le c\alpha^{n(x,y)}$. Then
\begin{equation}\label{eq:intbound}
\left|\int\phi\,d\mu'_{\eta,w}-\left( \int\phi\,d\mu+\delta\Big(S_{|w|}\phi(w)-|w|\int\phi\,d\mu\Big)\right)\right|
\le \delta c(2L+\tfrac 4{1-\alpha}),
\end{equation}
where $\mu'_{\eta,w}$ is given by
\eqref{eq:CoupleSplice}, $\delta=\nu_{\eta,2L+|w|+1}([1])$ as in the previous lemmas
and $S_{|w|}\phi(w)$ denotes $\inf_{x\in[w]}S_{|w|}\phi(x)$.
\end{lem}

\begin{proof}
We let $\bar X$ and $\bar\mu$ be as in the proof of Lemma \ref{lem:ent_estimate}.
As before, let $E=\{(x,y,z)\in\bar X\colon y_0=1\}$ and $r_E(x,y,z)$ denote the return time to $E$.
By definition, $\bar\mu(E)=\delta$, and $\delta$ is computed in \eqref{eq:delta}.
We let $\bar\mu_E$ be the normalized induced measure on the set $E$.
Notice that we have
\begin{equation}
\label{eq:coords}
\begin{split}
\int\phi\,d\mu&=\delta\int_E S_{r_E(\bar x)}\phi\circ\pi_1(\bar x)\,d\bar\mu_E(\bar x)\text{;\quad and}\\
\int\phi\,d\mu'_{\eta,w}&=\delta\int_E S_{r_E(\bar x)}\phi\circ\pi_3(\bar x)\,d\bar\mu_E(\bar x),
\end{split}
\end{equation}
where $\pi_1(x,y,z)=x$ and $\pi_3(x,y,z)=z$.

Let $m=|w|$. For $0\le j<L$ and $L+m\le j<2L+m$, we have
\begin{equation}\label{eq:transblock}
|\phi(T^j\pi_3(\bar x))-\phi(T^j\pi_1(\bar x))|\le c.
\end{equation}
If $\bar x=(x,y,z)\in E$, then $x_L\ldots x_{L+m-1}=w$.
If $x'_L\ldots x'_{L+m-1}=w$ also, we use the estimate $|\phi(T^ix)-\phi(T^ix')|\le 
c\alpha^{\min(|i-L+1|,|i-(L+m)|)}$,
valid for $i\le L\le L+m-1$ to obtain
$$
\left|\int_E \sum_{j=L}^{L+m-1}\phi(T^j \pi_3(\bar x))\,d\bar \mu_E-S_m\phi(w)\right|\le \frac{2c}{1-\alpha}.
$$
We also have
\begin{align*}
&\int_E\sum_{j=L}^{L+m-1}\phi(T^j \pi_1(\bar x))\,d\bar\mu_E=
\int_E S_m\phi(T^L(\pi_1(\bar x)))\,d\bar\mu_E\\
&=\frac1{\bar\mu(E)}\int_{\{\bar x\colon y_0=1\}}S_m\phi(T^Lx)\,d\bar\mu(x,y,z)\\
&=\int S_m\phi(T^Lx)\,d\mu(x)\\
&=m\int\phi(x)\,d\mu(x),
\end{align*}
where in the third line, we used the independence of $\pi_1(\bar x)$ and $\pi_2(\bar x)$.
Combining the two previous facts, we see
\begin{equation}\label{eq:Wblock}
\begin{split}
&\left|\int_E \big(S_m\phi(T^L\pi_3(\bar x))-S_m\phi(T^L\pi_1(\bar x))\big)\,d\bar\mu_E-
\left(
S_m\phi(w)-m\int\phi\,d\mu\right)\right|\\
&\le 2Lc+\frac{2c}{1-\alpha}.
\end{split}
\end{equation}
Finally, if $\bar x=(x,y,z)\in E$ and $2L+m\le j<r_E(\bar x)$, then
$$
|\phi(T^jz)-\phi(T^jx)|\le c \alpha^{\min(j-(2L+m-1),r_E(\bar x)-j)}.
$$
Summing the geometric series over this range of $j$'s yields
\begin{equation}\label{eq:commonparts}
\sum_{j=2L+m}^{r_E(\bar x)-1}
\Big|\phi(T^j\pi_3(\bar x))-
\phi(T^j\pi_1(\bar x))\Big|
\le \frac {2c}{1-\alpha}.
\end{equation}
Combining equations \eqref{eq:transblock}, \eqref{eq:Wblock} and \eqref{eq:commonparts}, we obtain
\begin{align*}
&\left|\int_E \big(S_{r_E(\bar x)}\phi\circ \pi_3 (\bar x)-S_{r_E(\bar x)}\phi\circ \pi_1 (\bar x)\big)\,d\bar\mu_E
-\left(S_m\phi(w)-m\int \phi\,d\mu\right)\right|
\\
&
 \le 2Lc+\frac {4c}{1-\alpha},
\end{align*}
so that the claimed result follows from \eqref{eq:coords}.
\end{proof}

\begin{lem}\label{lem:distort}
Let $\phi$ be a H\"older potential on a mixing shift of finite type $X$ and
let $\mu$ be the corresponding Gibbs measure with constant $C_\phi$
defined in \eqref{eq:Gibbs}.
For the partition $\cP$
and all $n\in\N$ we have
$$
nh_T(\mu)\le H_\mu(\mathcal P_n)\le nh_T(\mu)+\log C_\phi.
$$
\end{lem}

\begin{proof}
The lower bound is well known. For the upper bound, we have
$$
H(\mathcal P_n)=
-\int \log\mu(\mathcal P_n(x))\,d\mu(x),
$$
where $\mathcal P_n(x)$ is the element of $\mathcal P_n$ containing $x$.
By \eqref{eq:Gibbs},
$$
\left|
\log\mu(\mathcal P_n(x))-S_n\phi(x)+nP(\phi)\right|\le\log C_\phi,
$$
so that
$$
\left|-\int\log\mu(\mathcal P_n(x))\,d\mu(x)-\int\Big(nP_{\rm top}(\phi)-S_n\phi(x)\Big)\,d\mu\right|\le\log C_\phi,
$$
which yields
$$
\left|H_\mu(\mathcal P_n)-n\left(P_{\rm top}(\phi)-\int\phi\,d\mu\right)\right|\le \log C_\phi.
$$
Since $\mu$ is an equilibrium state for the potential $\phi$, we have the equality $P_{\rm top}(\phi)=h_T(\mu)+\int\phi\,d\mu$, so the above
gives $|H_\mu(\mathcal P_n)-nh_T(\mu)|\le\log C_\phi$ as claimed.
\end{proof}

The following
lemma controls the behaviour of $C_{t\phi}$ as $t$ runs over the reals.

\begin{lem}\label{lem:Ruelle}
Let $\phi$ be a H\"older continuous potential on a mixing shift of finite type $X$. 
Then there exist $a$ and $b$
such that $C_{t\phi}\le e^{a+b|t|}$ for all $t\in\R$.
\end{lem}

\begin{proof}
For the proof, we rely heavily on the results in Ruelle \cite[Chapter 5]{Ru1}.

We may assume without loss of generality that $\phi$ is one-sided. That is,
$\phi(x')=\phi(x)$ if $x'$ and $x$ satisfy $x'_i=x_i$ for all $i\ge 0$. 
To justify this reduction, by a result of Bowen \cite[Lemma 1.6]{Bo}, 
there is a H\"older continuous potential $\tilde\phi$ and a H\"older continuous
transfer function $h$ such that
$\tilde\phi=\phi+h-h\circ T$ and such that $\tilde\phi$ is one-sided.
Clearly, the equilibrium states for $t\tilde\phi$ and $t\phi$
agree for all $t\in\R$. Also $|S_n(t\phi)(x))-S_n(t\tilde\phi(x))|\le 2|t|\|h\|_\infty$,
so that from \eqref{eq:Gibbs}, we see 
$$
e^{-2|t|\|h\|_\infty}\le C_\phi/C_{\tilde\phi}\le
e^{2|t|\|h\|_\infty}.
$$

We now claim that $\phi$ may be expressed in the form $\sum_{n=0}^\infty \phi_n$,
where $\phi_n(x)$ depends only on $x_0,\ldots x_n$ and $\|\phi_n\|_\infty\le c
\alpha^n$ for each $n\ge 1$, where $c$ and $\alpha$ are the H\"older constant and exponent of $\phi$ respectively.
This proof is essentially contained in \cite{CoelhoQuas}. To establish the claim, let
$\psi_n(x)=\min\{\phi(y)\colon y_0...y_n = x_0...x_n\}$. Then for each $x$, sequence $(\psi_n(x))$ increases
to $\phi(x)$ and we check that $\psi_{n+1}(x)-\psi_n(x)\le \|\phi\|_\alpha\alpha^{n+1}$, 
where $\|\phi\|_\alpha=\max\{\|\phi\|_\infty,
\sup_{x\ne y}|\phi(x)-\phi(y)|/\alpha^{n(x,y)}\}$. 
Let $\phi_0=\psi_0$ and $\phi_n=\psi_n-\psi_{n-1}$ for $n\ge 1$. This gives the required 
decomposition of $\phi$. We abuse notation and write $\phi_n(x_0\ldots x_n)$
for the common value of $\phi_n(x)$ for all $x\in[x_0\ldots x_n]$. We write $X^+$
and $X^-$ for the projections of $X$ onto the positive and negative coordinates respectively. 

For comparison with Ruelle's book, the interaction $\Phi$ that he is considering is defined by
$$
\Phi_\Lambda(x_\Lambda)=\begin{cases}
\Phi_{[a,b]}(x_{[a,b]})=\phi_{b-a}(x_a\ldots x_b)&\text{if $\Lambda=[a,b]$};\\
0&\text{otherwise}
\end{cases}
$$
To see this, if we follow \cite[equation (3.3)]{Ru1} to convert $\Phi$ back to a potential,
the potential that we obtain is precisely $\phi$.

Following \cite[Sections 5.11--13]{Ru1}, we define a family of Perron-Frobenius operators,
mapping $C^\alpha(X^+)$ to itself by
$\mathcal P_tf(x)=\sum_{i\colon ix_0\in\mathcal \cL(X)} e^{t\phi(ix)}f(ix)$. 
By \cite[Proposition 5.13]{Ru1}, the leading eigenvalue of $\mathcal P_t$ is $P(t\phi)$.
We now apply \cite[Section 5.12]{Ru1} to obtain an expression for the
leading eigenfunction. Notice that the function $-tW_{\mathbb Z^-,\mathbb Z^{0+}}(\xi_<,\xi_\ge)$
(defined in \cite[equation (1.8)]{Ru1}) is equal to 
$t\sum_{\Lambda\supseteq\{-1,0\}}\Phi_\Lambda(\xi_<\xi_\ge |_\Lambda)$,
which may be expressed in our
notation as
$$
t\sum_{m=1}^\infty\sum_{n=0}^\infty \phi_{n+m}(y_{-m}\ldots
y_{-1}x_0\ldots x_n).
$$
The expression for the eigenfunction is then given by
\begin{equation}
\label{eq:efn}
\rho_t(x)=\frac 1K\int_{\mathsf{Pre}(x_0)}\exp\left(\sum_{m=1}^\infty\sum_{n=0}^\infty t\phi_{n+m}(y_{-m}\ldots
y_{-1}x_0\ldots x_n)\right)\,d\nu_t^-(y),
\end{equation}
where $K$ is a normalization constant, $\mathsf{Pre}(x_0)=\{y\in X^-\colon y_{-1}x_0\in \cL(X)\}$
and $\nu_t^-$ is the {\color{teal} (non-invariant)}``Gibbs state on $X^-$", given by 
$$
\nu_t^-([y_{-m}\ldots y_{-1}])=\lim_{n\to\infty}\nu_t^{(n)}([y_{-m}\ldots y_{-1}]),
$$
where
$$
\nu_t^{(n)}([y_{-m}\ldots y_{-1}])=
\frac
{
\sum_{zy\in \cL_n(X)}Q_t^{(n)}(z_{-n}\ldots z_{-m-1}y_{-m}\ldots y_{-1})
}
{
\sum_{z\in \cL_n(X)}Q_t^{(n)}(z_{-n}\ldots z_{-1})
},
$$
and
$$
Q_t^{(n)}(z_{-n}\ldots z_{-1})=
\exp\left(t\sum_{1\le i\le j\le n}\phi_{j-i}(z_{-j}\ldots z_{-i})\right).
$$
To see the convergence of the limit defining $\nu_t^-$, note that
$\nu^{(n)}_t$ is exactly the finite Gibbs state on $X_{[-n,-1]}$ for the interaction $t\Phi$.
By Theorem 5.3(b) of Ruelle, there is a unique Gibbs state on $X^-$ for the interaction $t\Phi$.
By Section 1.4 and 1.6, any subsequential limit of $(\nu_t^{(n)})_{n\ge 1}$
is a Gibbs state on $X^-$ for $t\Phi$. Therefore the finite Gibbs states converge in the weak$^*$-topology
to $\nu_t^-$. \\
If $L$ is the mixing length of the shift and $n>L+1$,
we claim that 
\begin{equation}\label{eq:Gibbsclaim}
\log\left(\frac{\nu_t^{(n)}([i])}{\nu_t^{(n)}([j])}\right)
\le L\log|\cA|+2t(L+1)\sum_{n=0}^\infty\|\phi_n\|.
\end{equation}
Assuming this inequality for now, by taking the limit as $n\to\infty$, we will obtain
\begin{equation}\label{eq:Gibbs-}
\left|\log\left(\frac{\nu_t^-([i])}{\nu_t^-([j])}\right)\right|\le a+b|t|,
\end{equation}
where $a$ and $b$ only depend on $\phi$ and $X$. 
An identical argument gives
\begin{equation}\label{eq:Gibbs+}
\left|\log\left(\frac{\nu_t^+([i])}{\nu_t^+([j])}\right)\right|\le a+b|t|,
\end{equation}
where $\nu_t^+$ is the (non-invariant) Gibbs state on $X^+$.

To show \eqref{eq:Gibbsclaim}, we let $n=m+L+1$ with $m>0$ and observe
\begin{equation*}\label{eq:Gibbsobs}
\frac{\nu_t^{(n)}([i])}{\nu_t^{(n)}([j])}=
\frac{\sum_{x\in\cL_m}\sum_{y\in \cL_L: xyi\in\cL}Q_t^{(n)}(xyi)}
{\sum_{x\in\cL_m}\sum_{y'\in \cL_L: xy'j\in\cL}Q_t^{(n)}(xy'j)}.
\end{equation*}

To bound this quantity, it suffices to obtain a bound for 
\begin{equation*}\label{eq:Gibbs_to_bound}
\frac{\sum_{y\in \cL_L: xyi\in\cL}Q_t^{(n)}(xyi)}
{\sum_{y'\in \cL_L: xy'j\in\cL}Q_t^{(n)}(xy'j)},
\end{equation*}
that does not depend on $x\in\cL_m$.
We bound this by estimating the ratio $Q_t^{(n)}(xyi)/Q_t^{(n)}(xy'j)$
for any $y$ and $y'$ satisfying the conditions.
From the definition
$$
\frac{Q_t^{(n)}(xyi)}{Q_t^{(n)}(xy'j)}=\exp\left(t\sum_{1\le p\le q\le n}\Big(\phi_{q-p}(xyi_{[p,q]})-\phi_{q-p}(xy'j)_{[p,q]}\Big)\right).
$$
Clearly the terms involving only $x$'s cancel so we obtain
$$
\frac{Q_t^{(n)}(xyi)}{Q_t^{(n)}(xy'j)}\le \exp\left(2t(L+1)\sum_{k=0}^\infty \|\phi_k\|\right),
$$
since there are at most $L+1$ terms in the sum involving $\phi_k$ that do not involve only $x$ terms.
We therefore see that for each $x\in \cL_m$,
\begin{equation*}
\frac{\sum_{y\in \cL_L: xyi\in\cL}Q_t^{(n)}(xyi)}
{\sum_{y'\in \cL_L: xy'j\in\cL}Q_t^{(n)}(xy'j)}
\le |\cA|^L \exp\left(2t(L+1)\sum_{k=0}^\infty \|\phi_k\|\right),
\end{equation*}
and we conclude that \eqref{eq:Gibbsclaim} and \eqref{eq:Gibbs-} hold. 

Looking at \eqref{eq:efn}, we define $A$ by
\begin{align*}
A=\left|\sum_{m=1}^\infty\sum_{n=0}^\infty \phi_{n+m}(y_{-m}\ldots
y_{-1}x_0\ldots x_n)\right|&\le \sum_{m=1}^\infty\sum_{n=0}^\infty \|\phi_{n+m}\|\\
&=\sum_{k=1}^\infty k\|\phi_k\|.
\end{align*}
Then the integrand in \eqref{eq:efn} takes values in the range $[e^{-|t|A},e^{|t|A}]$,
so that $\rho_t(x)$ is in the range $[\nu_t^-(\mathsf{Pre}(x_0))e^{-t|A|},[\nu_t^-(\mathsf{Pre}(x_0))e^{t|A|}]$.
Using \eqref{eq:Gibbs-}, it follows that there are constants $a$ and $b$ such that
\begin{equation}\label{eq:rhoratio}
\left|\log\left(\frac{\rho_t(x)}{\rho_t(x')}\right)\right|\le a+b|t|.
\end{equation}

As mentioned in \cite[Section 5.14]{Ru1},
$\chi_t=t\phi+\log\rho_t-\log \rho_t\circ T-P(t\phi)$ is a normalized potential.
That is, it satisfies $\sum_{ix_0\in\cL(X)}\exp\chi_t(ix)=1$ for all $x\in X^+$.
Following Ruelle, we define the Perron-Frobenius operator of $\chi_t$ by 
$$
\mathcal S_tf(x)=\sum_{i\in\cA, ix\in\cL}e^{\chi_t(ix)}f(ix).
$$
By the normalization condition, $\mathcal S_t\mathbf 1=\mathbf 1$ (see Theorem 5.14).
One can verify from \cite[(5.18)]{Ru1} that $\mathcal S_t^*\mu_t^+=\mu_t^+$, where 
$\mu_t^+$ is the one-sided version of $\mu_t$. Combining this with \cite[(5.22)]{Ru1}
and letting $C=[u_0\ldots u_{n-1}]$, we obtain
\begin{align*}
\mu_t(C)&=\int \mathbf 1_C\,d\mu_t\\
&=\int \mathcal S_t^n(\mathbf 1_C)\,d\mu_t\\
&=\int_{\mathsf{Suc}(u_{n-1})}\exp S_n\chi_t(ux)\,d\mu_t(x),
\end{align*}
where $\mathsf{Suc}(u_{n-1})=\{x\in X^+\colon u_{n-1}x\in X^+\}$.
It follows that
$$
\log \mu_t([x_0\ldots x_{n-1}])\le S_n\chi_t(x)+\operatorname{var}_n S_n\chi_t(x)
$$
and
$$
\log \mu_t([x_0\ldots x_{n-1}])\ge S_n\chi_t(x)-\operatorname{var}_n S_n\chi_t(x)
+\min_i\log\mu_t([i]).
$$
By \cite[(5.18)]{Ru1}, 
$$
\mu_t([i])=\int_{[i]}\rho_t(x)\,d\nu_t^+(x).
$$
Using \eqref{eq:Gibbs+} and \eqref{eq:rhoratio}, we see there exists $a,b$ such that $|\log\mu_t([i])|\le a+b|t|$
for each $i$ and $t$.

From the cohomological relationship, we see
$\operatorname{var}_n S_n\chi_t\le t\operatorname{var}_n S_n\phi + 2\|\log\rho_t\|$.
We can compute $\operatorname{var}_n S_n\phi\le 
\sum_{i=0}^{n-1}\sum_{j=n-i}^\infty\|\phi_j\|_\infty$, so that $\operatorname{var}_n S_n\phi$
is uniformly bounded in $n$. 
It follows that there exist $c$ and $d$, depending only on $X$ and $\phi$ such that
$$
\Big|\log \mu_t([x_0\ldots x_{n-1}])-(tS_n\phi(x)-n\log\lambda_t)\Big|\le c+d|t|,
$$
as required.
\end{proof}
\color{black}

\section{Proofs of Main theorems}\label{sec:main_proofs}
We now turn to the proofs of the main theorems, which mostly consist of
estimating quantities of the form $h_T(\mu')+s\int\phi\,d\mu'$ from below where $\mu'$ is one of the measures built
in the previous section (and $\mu$ is $\mu_\infty$ in the case of Theorem \ref{thm:gap_at_infinity}
or $\mu_{t}$ in the case of Theorem \ref{thm:strictconv}). This then gives a lower
bound for $p_\phi(s)$. To bound $h_T(\mu')+s\int\phi\,d\mu'$
from below, we rely on bounds from the previous section,
showing that
$h_T(\mu')$ exceeds $h_T(\mu)$ plus a term of order $-\eta\log\eta$
and that $s\int\phi\,d\mu'$ is at least $s\int \phi\,d\mu$ minus a term of order $s\eta$).
If $\eta$ is taken to be less than $e^{-Cs}$ for a suitable $C$,
the gain dominates the loss by an amount of order $\eta$.

We restate Theorem \ref{thm:gap_at_infinity} for convenience.

\begin{thm*}
Let $(X,T)$ be a mixing subshift of finite type with positive entropy.
Let $\phi$ be a H\"older potential
that is not cohomologous to a constant.
Then there exist $C>0$ and $t_0$ such that $p_\phi(t)\ge\ell_\infty(t)+e^{-Ct}$
for all $t\ge t_0$, where $\ell_\infty(t)$ is the
asymptote to $p_\phi$ at infinity.
\end{thm*}

\begin{proof}
By \cite{Bousch-bilateral}, there exists a H\"older continuous function
$\psi$ that is cohomologous to $\phi$ such that
$A(\phi) \le \psi(x) \le B(\phi)$ for
all $x\in X$, where
$A(\phi)=\min_{\nu\in \invmeas}\int \phi\,d\nu$ and
$B(\phi)=\max_{\nu\in\invmeas}\int \phi\,d\nu$. The assumption that
$\phi$ is not cohomologous to a constant implies that $A(\phi)\ne B(\phi)$. Since the pressure functions of $\phi$ and $\psi$ coincide, we derive the required estimate for $p_\psi(t)$.

Let $\mu_\infty$ be a measure achieving the slant asymptote $\ell_\infty$, so that the support of
$\mu_\infty$ is contained in the proper  subset of $X$:
$$
S_\text{max}:=\{x\in X\colon \psi(x)=B(\phi)\}.
$$
As before, denote by $L$ the mixing length of $X$. Let $[w]$ be a cylinder set lying in the complement of $\text{supp}(\mu_\infty)$ where $w$
is a word of length at least $3L$ with no long overlaps. Such a word always exists, see e.g. \cite[Theorem 8.3.9]{Lothaire2}.
Let $M=|w|+2L+1$. The word $w$ is now fixed for the remainder of this proof.
For any $\eta>0$, we equip $Y=\{0,1\}^\Z$ with the measure $\nu=\nu_{\eta,M}$ 
constructed in the previous section and use
$\nu$ to build a measure $\mu'_{\eta,w}$ as in \eqref{eq:CoupleSplice} (where $\mu$ is taken to be
the measure $\mu_\infty$).

We shall show that there exist constants $a$ and $b$ (depending only on $|w|$, $L$, and the H\"older constants of $\psi$, ) 
such that for all sufficiently large $t$, there is an $\eta>0$ such that
$$
h_T(\mu'_{\eta,w})+t\int \psi\,d\mu'_{\eta,w}\ge \ell_\infty(t)+ae^{-bt},
$$
For now, we leave $\eta\in(0,1)$ arbitrary, and choose it later in the proof. 
Let $\delta=\nu_{\eta,M}([1])$, so that $\delta=(M+(e^\eta-1)^{-1})^{-1}$
by Lemma \ref{lem:build-nu}.

We see from Lemma \ref{lem:int_estimate} that
\begin{equation*}
\int\psi\,d\mu'_{\eta,w} \ge \int \psi\,d\mu_\infty-\delta c\left(M-1+\frac{4}{1-\alpha}\right),
\end{equation*}
where we recall that $M-1=2L+|w|$, and where $c$ and $\alpha$ are the constants in the H\"{o}lder condition for $\psi$.

Since $\mu_\infty([w])=0$ the condition $\delta\ge\frac 15\mu([w])$ holds for all $\eta>0$.
 The estimate in Lemma \ref{lem:ent_estimate}  works whenever $\eta\in [\mu([w]),1)$, and since here $\mu([w])=0$, it works for all $\eta\in(0,1)$
Combining it with the entropy formula for the
measure $\nu_{\eta,M}$ from Lemma \ref{lem:build-nu}, we conclude that 
\begin{align*}
h_T(\mu'_{\eta,w})&\ge h_T(\mu_\infty)+\delta\big(-\log\eta-H_{\mu_\infty}(\mathcal P_{M-1})-8\log 2\big)\\
&\ge h_T(\mu_\infty)+\delta\big(-\log\eta-((M-1)\log|\cA|+8\log 2)\big).
\end{align*}
It follows that
\begin{align*}
P_{\rm top}(t\psi)&\ge h_T(\mu'_{\eta,w})+t\int\psi\,d\mu'_{\eta,w}\\
&\ge \left(h_T(\mu_\infty)+t\int\phi\,d\mu_\infty\right)+\delta\big(-\log\eta-(c_1+c_2|t|)\big),
\end{align*}
where $c_1=(M-1)\log|\cA|+8\log 2$ and $c_2=c(M-1+\frac 4{1-\alpha})$. 
This bound holds for all $\eta\in(0,1)$.
We have $\delta\ge \frac\eta 2$ for small values of $\eta$ by \eqref{eq:delta}.
Substituting $\eta=\exp(-1-c_1-c_2|t|)$, we
deduce that $p_\phi(t)=p_\psi(t)\ge \ell_\infty(t)+\frac 12e^{-(c_1+1)}e^{-c_2|t|}$ 
for all large values of $t$, as required.
\end{proof}

For locally constant functions, the
true gap between $p_\phi(t)$ and $\ell_\infty(t)$ is asymptotically exponential, matching the form
of the lower bound in the previous theorem:
\begin{example}\label{ex:optimal_bnd}
Let $(X,T)$ be the full two-sided shift on the alphabet $\{1,...,k\}$ and $\phi:X\to\R$ be a potential which is constant on cylinders of length 1, i.e. $\phi(x)=c_{x_0}$ where $c_1,\ldots,c_k$ are fixed real numbers. Then
$p_\phi(t)=\log(e^{c_1t}+\ldots+e^{c_{k}t})$.

In the case $c_i>\max_{j\ne i}c_j$, $p_\phi(t)=c_it+O(e^{-\Delta t} )$ where $\Delta=c_i-\max_{j\ne i}c_j$.
\end{example}

We now restate Theorem \ref{thm:strictconv}.
\begin{thm*}
Let $(X,T)$ be a mixing subshift of finite type with positive entropy.
Let $\phi$ be a H\"older potential
that is not cohomologous to a constant. Then there exists $c>0$
such that for any $t\in\R$, and any $h\in (0,1)$,
\begin{equation*}
p_\phi(t+h)+p_\phi(t-h)-2p_\phi(t)>e^{-c(1+t^2)/h}.
\end{equation*}
\end{thm*}

\begin{proof}
We write $\ell_{t}(s)=h(\mu_{t})+s\int \phi\,d\mu_{t}$, where $\mu_t$ denotes the (unique) equilibrium state for
the H\"older continuous potential $t\phi$. This is the tangent line to $p_\phi$ at $(t,p_\phi(t))$. We also define
the ``gap function", $g_t(s)=p_\phi(s)-\ell_t(s)$. Since $p_\phi$ is known to be strictly convex, we see $g_t(t)=0$ and 
$g_t(s)>0$ for $s\ne t$. 
Since $\ell_t$ is affine, we see
\begin{align*}
p_\phi(t+h)+p_\phi(t-h)-2p_\phi(t)&=g_t(t+h)+g_t(t-h)-2g_t(t)\\
&=g_t(t+h)+g_t(t-h).
\end{align*}
Since $g_t(s)\ge 0$ for each $s$, to demonstrate the theorem it suffices 
to show that there exists $c>0$ such that
for all $t\in\R$ and all $h\in (0,1)$, 
\begin{equation}\label{eq:toshow}
\max(g_t(t+h),g_t(t-h))>e^{-c(1+t^2)/h}.
\end{equation}
As before, let $A=\min_\mu\int\phi\,d\mu$ and $B=\max_\mu\int\phi\,d\mu$, where $\mu$ runs over the collection of
invariant measures. 
By \cite{Ru1}, $\int\phi\,d\mu_t=p_\phi'(t)$ for all $t\in\R$ and by convexity
of pressure, $t\to\int\phi\,d\mu_t$ is an increasing function.
Recall from the folklore theorem on pressure asymptotics that
$p_\phi'(t)\to A$ as $t\to-\infty$ and $p_\phi'(t)\to B$ as $t\to\infty$. 
Let $C=\int\phi\,d\mu_0$ and let $\gamma=\min(\frac 12(B-C),\frac 12(C-A))$. 

As mentioned in the proof of Theorem \ref{thm:gap_at_infinity} there exists a H\"older continuous
function $f$ such that $\psi=\phi+f\circ T-f$ takes values in the range $[A,B]$.
It follows that there exist (disjoint) subshifts $X_\alpha$ and $X_\beta$ on which $\psi$ takes the constant
values $A$ and $B$ respectively. Since the pressure functions of $\phi$ and $\psi$
coincide, we notice $g_t(s)=p_\psi(s)-\ell_t(s)$ for all $s$. 
Moreover, $\mu_t$ is the unique equilibrium state for $t\psi$ for all $t$. We now use $\psi$ to estimate the gap function.

A word $w$ of length $m$ is called \emph{heavy} if $S_m\psi(x)>m(B-\gamma)$ for all $x\in[w]$
and similarly a word $w$ is \emph{light} if $S_m\psi(x)<m(A+\gamma )$ for all $x\in [w]$.
All sufficiently long words in the languages of $X_\alpha$ or $X_\beta$ are light and heavy respectively
by continuity of $\psi$. Indeed, the same is true for long words that agree with words in 
$X_\alpha$ or $X_\beta$ up to a small number of modifications.

In particular, one may find a pair of heavy words $u_h$ and $v_h$ of the same length
such that arbitrary concatenations
of $u_h$ and $v_h$ are legal in $X$ and so that $u_h$ does not appear as a sub-word of the infinite
concatenation of $v_h$'s and $v_h$ does not appear as a sub-word of the infinite concatenation of $u_h$'s.
Likewise there exist light words $u_l$ and $v_l$ with the analogous properties.
The four words $u_l,u_h,v_l,v_h$ are now fixed for the remainder of the proof.

The proof is divided into two cases, according to whether $t\ge 0$ or $t<0$.
If $t<0$, then $\int\psi\,d\mu_t<C$ and, similarly 
to the proof of Theorem \ref{thm:gap_at_infinity}, we 
overwrite random parts of $\mu_t$-typical points, inserting additional heavy words, and obtain a lower 
bound for $g_t(t+h)$. In the other case, we insert additional light words in $\mu_t$-typical points, 
obtaining a lower bound for $g_t(t-h)$.

Suppose for now that $t<0$, so that $\int\psi\,d\mu_t<C$. We write $\mu$ for $\mu_t$ for the remainder of the proof. 
For each $k\ge 2$ such that $2k|u_h|\ge L$, we let $w_k=u_h^kv_h^k$ and 
define a sequence of measures $\mu'_k$ using the coupling and splicing construction
\eqref{eq:CoupleSplice} by $\mu'_k=\mu'_{\eta_k,w_k}$, where $\eta_k=\mu([w_k])$.
By the Variational Principle, we have 
\begin{align*}
g_t(t+h)&=P_{\rm top}((t+h)\psi)-\ell_t(t+h)\\
&\ge \sup_k\left(h_T(\mu'_k)+(t+h)\int \psi\,d\mu'_k-h_T(\mu)-(t+h)\int \psi\,d\mu\right)
\end{align*}
The strategy of the proof is to look for a suitable $k$ to obtain the desired lower bound.

As before, let $\delta_k=\nu_{\eta_k,w_k}([1])$.
We then use the previous lemmas to estimate the gap
$(h_T(\mu_k')-h_T(\mu))+(t+h)(\int\psi\,d\mu_k'-\int\psi\,d\mu)$.

From Lemmas \ref{lem:build-nu} and \ref{lem:ent_estimate}, 
\begin{equation}\label{eq:hest1}
h_T(\mu_k')-h_T(\mu)\ge \delta_k(-\log\eta_k-H_\mu(\mathcal P_{2L+|w_k|})-8\log 2).
\end{equation}
Using Lemmas \ref{lem:distort} and \ref{lem:Ruelle}, 
\begin{equation}\label{eq:hest2}
H_\mu(\mathcal P_{2L+|w_k|})\le (2L+|w_k|)h_T(\mu)+a+b|t|,
\end{equation}
where the constants $a$ and $b$ are from Lemma \ref{lem:Ruelle}.
Using the Gibbs inequality, Lemma \ref{lem:Ruelle} and the Variational Principle, we obtain 
$$
\eta_k=\mu([w_k])\le e^{a+b|t|}e^{S_{|w_k|}(t\psi)(w_k)-|w_k|(h_T(\mu)+\int(t\psi)\,d\mu)}.
$$
Hence 
\begin{equation}\label{eq:hest3}
-\log\eta_k\ge -t\left(S_{|w_k|}\psi(w_k)-|w_k|\int\psi\,d\mu\right)+|w_k|h_T(\mu)-a-b|t|.
\end{equation}
Combining \eqref{eq:hest1}, \eqref{eq:hest2} and \eqref{eq:hest3}, we see
\begin{equation}\label{eq:hdiff}
\begin{split}
&h_T(\mu_k')-h_T(\mu)\\
&\ge \delta_k\left(-t\left(S_{|w_k|}\psi(w_k)-|w_k|\int\psi\,d\mu\right)
-2Lh_T(\mu)-2a-2b|t|-8\log 2\right).
\end{split}
\end{equation}
From Lemma \ref{lem:int_estimate},
\begin{equation}\label{eq:intdiff}
\begin{split}
(t+h)\left(\int\psi\,d\mu_k'-\int\psi\,d\mu\right)
&\ge
\delta_k\left((t+h)\left(S_{|w_k|}\psi(w_k)-|w_k|\int\psi\,d\mu\right)\right.\\
&\left.-2c(|t|+1)\left(L+\frac 2{1-\alpha}\right)\right),
\end{split}
\end{equation}
where $c$ and $\alpha$ are as in the statement of the lemma.

Adding \eqref{eq:hdiff} and \eqref{eq:intdiff}, we see 
\begin{equation}\label{eq:gapest}
g_t(t+h)\ge \delta_k \left(h\left(S_{|w_k|}\psi(w_k)-|w_k|\int\psi\,d\mu\right)-(a'+b'|t|)\right),
\end{equation}
where $a'=2a+8\log 2 +2Lh_T(\mu)+2c(L+\frac 2{1-\alpha})$ and 
$b'=2b +2c(L+\frac 2{1-\alpha})$. Notice that $a'$ and $b'$ depend on
$\psi$ and $X$, but not on $t$ or $h$. 
By the heaviness condition and the assumption that $\int\psi\,d\mu<C$, 
we see that $S_{|w_k|}\psi(w_k)-|w_k|\int\psi\,d\mu> |w_k|((B-\gamma)-C)\ge \gamma|w_k|$, so we have shown
$$
g_t(t+h)\ge \delta_k(ht\gamma |w_k|-(a'+b'|t|)),
$$ 
where this inequality holds for each $k$.
Recalling that $|w_k|=2k|u_h|$ and taking $k=\lceil (1+a'+b'|t|)/(2\gamma h|u_h|)\rceil$, we obtain 
$g_t(t+h)\ge \delta_k$. Since we chose $\eta_k=\mu([w_k])$, 
Lemma \ref{lem:ent_estimate} implies $\delta_k\ge \frac 15\mu([w_k])$.
We then use the Gibbs inequality one more time to estimate $\mu([w_k])$ from below.
Taking $x\in [w_k]$, we get
\begin{align*}
    \log\mu([w_k])&\ge (-a-b|t|)+S_{|w_k|}t\psi(x)-|w_k|P_{\rm top}(t\psi)\\
    &=(-a-b|t|)+t\left(S_{|w_k|}\psi(x)-|w_k|\int\psi\,d\mu\right)-|w_k|h_\sigma(\mu)\\
    &\ge (-a-b|t|)-\gamma |t| |w_k|-|w_k|h_\sigma(\mu),\end{align*}
where we used $t<0$, $S_{|w_k|}\psi(x)\ge (C+\gamma)|w_k|$ (by heaviness of $w_k$)
and $\int\psi\,d\mu<C$.
Since $|w_k|\le (a''+b''|t|)/h$, where $a''=2|u_k|+(1+a')/\gamma$ and $b''=b'/\gamma$,
we obtain
$$
\log\mu([w_k])\ge -a-b|t|-(\log|\cA|+\gamma |t|)(a''+b''|t|)/h).
$$
Expanding, using $0<h<1$ and $|t|<\frac 12(1+t^2)$, we see there exists $c>0$ such that
$$
\mu([w_k])\ge e^{-c(1+t^2)/h}.
$$


In the case where $t\ge 0$, as mentioned above, we obtain a lower bound 
for $g_t(t-h)$ by inserting additional light words $w_k=u_l^kv_l^k$.
The inequality \eqref{eq:hdiff} is unchanged; the inequality \eqref{eq:intdiff}
becomes
\begin{equation}\label{eq:intdiff2}
\begin{split}
(t-h)\left(\int\psi\,d\mu_k'-\int\psi\,d\mu\right)&\ge
\delta_k\left((t-h)\left(S_{|w_k|}\psi(w_k)-|w_k|\int\psi\,d\mu\right)\right.\\
&\left.-2c(|t|+1)\left(L+\frac 2{1-\alpha}\right)\right),
\end{split}
\end{equation}
Combining \eqref{eq:hdiff} and \eqref{eq:intdiff2} gives
\begin{equation}\label{eq:gapest2}
g_t(t-h)\ge \delta_k \left(-h\left(S_{|w_k|}\psi(w_k)-|w_k|\int\psi\,d\mu\right)-(a'+b'|t|)\right),
\end{equation}
By lightness of $u_l$ and $v_l$,  $S_{|w_k|}\psi(w_k)-|w_k|\int\psi\,d\mu<-\gamma|w_k|$. We then
obtain a similar inequality to the 
earlier one:
$g_t(t-h)\ge \delta_k(h\gamma|w_k|-(a'+b'|t|))$ and the remainder of the proof is as before.
\end{proof}
\color{black}

\section{Arbitrarily slow convergence}\label{sec:slow_conv}

So far, we have proved that the pressure function cannot approach its asymptote ``too fast"
by providing an exponential lower bound on the gap. In this section we show the
absence of a corresponding upper bound. For a given mixing subshift of
finite type we construct potentials for which the convergence of the pressure to the
asymptote is arbitrarily slow, thus proving Theorem \ref{thm:large_gap}.

First, we construct a subshift $Y$ of a given entropy where a fixed word $u$ appears periodically.
For that, we make use of $\beta$-shifts. For $\beta>1$ the  $\beta$-shift is defined as the 
smallest two-sided subshift of $\{0,...,\lceil\beta\rceil-1\}^{\Z}$ which contains all 
sequences of the coefficients in $\beta$-expansions of real numbers in $[0,1)$, see e.g. \cite{Parry}. 
It is well known that the entropy of the $\beta$-shift is $\log \beta$ and there is a unique measure 
of maximal entropy which is fully supported.


\begin{lem}\label{lem:Yconstruction}
Let $(X,T)$ be a mixing subshift of finite type with positive entropy, let $0\le b< h_{\rm top}(X)$
and let $u\in \cL(X)$.
There exists an $\ell>0$ and a subshift $Y$ with the following properties:
\begin{enumerate}
\item The alphabet of $Y$ is a subset of $\cL_\ell(X)$,
where each symbol $w\in \cA(Y)$ is a word in $\cL_\ell(X)$ of the form $uv$
with the additional property that $vu\in\cL_\ell(X)$; 
\item $h_\text{top}(Y)=b\ell$.
\end{enumerate}
\end{lem}

The first condition ensures that the symbols appearing in $Y$ may be
concatenated to give a point of $X$ (with $u$'s occurring every $\ell$ steps); the
second condition will ensure that the subshift of $X$ formed from
these concatenations has topological entropy $b$.

%

\begin{proof}
Since $b<h_\text{top}(X)$, and $\#\cL_k(X)\ge e^{k\,h_\text{top}(X)}$ for all $k$,
there exists an $n$ such that $\#\cL_n(X)>e^{b(n+|u|+2L)}$ where
$L$ is the mixing length of $X$. Let $\ell=n+|u|+2L$ and set $N=\lceil e^{\ell b}\rceil$
(so that $N\le \#\cL_n(X)$).
Enumerate a subcollection of $N$ elements of $\cL_n(X)$ as $v_0,\ldots,v_{N-1}$.
By the mixing condition, there exist words
$p_1,\ldots,p_N$ and $q_1,\ldots,q_N$ in $\cL_L(X)$ such that $w_i:=up_iv_iq_i\in \cL_\ell(X)$ and
$q_iu\in\cL_{L+|u|}(X)$ for each $i$. Let $\cW=\{w_0,\ldots,w_{N-1}\}$.
By construction, these words are distinct.
Let $B$ denote the standard $\beta$-shift on the alphabet $\{0,\ldots,N-1\}$ where $\beta=e^{b\ell}$.
As was mentioned above, $h_\text{top}(B)=\log\beta=b\ell$.
We then let $Y$ be the image of $B$ under the bijective
one-block map $\theta\colon B\to\cW^\Z$ defined by $\theta(b)_j=w_{b_j}$, so that
$Y$ satisfies the conditions in the statement of the lemma.
\end{proof}

In the next lemma we describe an auxiliary shift $Z_\infty\subset\{0,1\}^\Z$ required for
our construction. A shift of this type was used by Rothstein in \cite{Rothstein} to produce
an example of a loosely Bernoulli process which is not very weak Bernoulli.
Hence, we will refer to $Z_\infty$ as the Rothstein shift.

\begin{lem}\label{lem:Rothstein}
Given an increasing sequence of positive integers $(n_j)_{j=1}^\infty$ there is a sequence
of nested subshifts $Z_j$ of $\left(\{0,1\}^\Z,\sigma\right)$ and a shift $Z_\infty=\bigcap Z_j$ satisfying
  \begin{itemize}
    \item $\displaystyle{h_{\rm top}(Z_j)=\frac{\log 2}{2^j}}$ (so that $h_{\rm top}(Z_\infty)=0$);
    \item For $j\in\N$ and $z\in Z_j$ we have $\displaystyle {d(z,Z_\infty)<2^{-n_{j+1}}}$.
  \end{itemize}
\end{lem}

\begin{proof}
  For a set of words $W$ on alphabet $\{0,1\}$ denote by $W^n$ the set of all possible concatenations of $n$ elements of $W$. We define a sequence $(W_j)$ of sets of words inductively. Let $(n_j)_{j=1}^\infty$ be an increasing sequence of positive integers. We start with $W_0=\{0,1\}$, which are two words of length 1. We set
$$W_{j+1}=\big\{ww:w\in W_j^{n_{j+1}}\big\}\text{ for }j=0,1,2...$$
Denote by $l_j$ the length of words in $W_j$. Then $l_0=1$ and one can easily check that $l_j=2^jn_1...n_j$. Also, since $W_0$ consists of two words, we see that $|W_1|$ has $2^{n_1}$ words, and the number of words in $W_j$ is
$$|W_j|=|W_{j-1}|^{n_j}=|W_{j-2}|^{n_{j-1}n_j}=...=2^{n_1n_2...n_j}.$$

Let $Z_j$ be the subshift which consists of all possible concatenations of words from $W_j$. Note that $_j$ is a sofic shift, since it is presented by a finite labeled graph with a loop corresponding to each word in $W_j$ (see e.g. \cite[Chapter 4]{LindMarcus}). Hence, $Z_j$ has a unique measure of maximal entropy. Since there are $|W_j|^n$ words of length $nl_j$ in $Z_j$, we compute
\begin{equation}\label{eq:RothEntropy}
  h_{\rm top}(Z_j)=\lim_{n\to\infty}\frac{\log (l_j2^{n(n_1...n_j)})}{nl_j}=\frac{n_1...n_j\log 2}{l_j}=\frac{n_1...n_j\log 2}{2^j n_1...n_j}=\frac{\log 2}{2^j}.
\end{equation}
The Rothstein shift $Z_\infty$ is defined as the set of all elements $x\in\{0,1\}^\Z$
such that every finite word in $x$ is a subword of $w\in W_j$ for some $j$. It follows from (\ref{eq:RothEntropy}) that $h_{\rm top}(Z_\infty)=0$.

We need an estimate on the distance between points in $Z_j$ and the shift $Z_\infty$.
Suppose $z\in Z_j$ for some $j\in\N$. We claim that for $n=\frac12(l_{j+1}-l_j)$ the
word $z_{-n}\ldots z_{n-1}$ is a subword of a word in $W_{j+2}$.
Indeed, $z_{-n}\ldots z_{n-1}$ has length $(2n_{j+1}-1)l_j$ and is in the
language of $Z_j$. Therefore, there are words $w_{i}\in W_j$ with $i=-n_{j+1},...,n_{j+1}-1$
such that $z_{-n}...z_{n-1}$ is a subword of their concatenation $w_{-n_{j+1}}\cdots w_{n_{j+1}-1}$. Denote $u=w_{-n_{j+1}}\cdots w_{-1}$ and $v=w_{0}\cdots w_{n_{j+1}-1}$, so that $z_{-n}...z_{n-1}$ is a subword of $uv$. Then by definition of $W_{j+1}$, $uu,vv\in W_{j+1}$ and hence $uuvv=(uu)(vv)$
is a subword of a word in $W_{j+2}$ as a concatenation of two words in $W_{j+1}$.
It follows that $z_{-n}\ldots z_{n-1}$ is a subword of a word in $W_{j+2}$
and hence a subword of a word in $W_k$ for all $k\ge j+2$. It follows that
$z_{-n}\ldots z_{n-1}\in \cL(Z_\infty)$. We deduce
\begin{equation}\label{eq:RothDistance}
  d(z,Z_\infty)\le 2^{-\frac{l_{j+1}-l_j}{2}},\text{ whenever } z\in Z_j.
\end{equation}
Recall that $l_j=2^jn_1..n_j$, so that $\frac12(l_{j+1}-l_j)=2^{j-1}n_1...n_j(2n_{j+1}-1)\ge n_{j+1}$. It follows from above that for $z\in Z_j$ we have $d(z,Z)\le 2^{-n_{j+1}}$, as required.
\end{proof}

We are now ready to provide the proof of Theorem \ref{thm:large_gap}. We remind the reader
of the statement.
\begin{thm*}
Suppose $(X,T)$ is any mixing subshift of finite type with positive entropy and
$f:\R\to\R$ is any convex function asymptotic to a line $\ell_\infty(t):=at+b$,
where $0\le b< h_{\rm top}(X)$.  Then there exists a H\"{o}lder potential $\phi:X\to \R$ such that
\begin{itemize}
  \item[(1)]  $p_\phi(t)$ is asymptotic to the line $\ell_\infty(t)=at+b$ as $t\to\infty$;
  \item[(2)] $p_{\phi}(t)>f(t)$ for all sufficiently large $t$.
 \end{itemize}
\end{thm*}

Before beginning the proof, we show why the result holds in a special case,
the case where $b=0$ and $X$ is a full shift on 2 symbols. 
Let $f(t)$ be a positive function decreasing to 0. 
We construct a
sequence $(n_j)$ such that if $Z_\infty$ is the associated Rothstein shift and
$\phi(x)=-d(x,Z_\infty)$, then $P(t\phi)\ge f(t)$ for all sufficiently large $t$. 

Let $j_0$ be such that $f(1)> 1/2^{j_0+1}$. Let $n_j=1$ for $j\le j_0$. 
For each $j\ge j_0$, let $t_j$ be such that $f(t_j)=1/2^{j+1}$. 
For each $j\ge j_0$, choose $n_{j+1}$ such that
$$
\frac{\log 2}{2^j}-\frac{t_{j+1}}{2^{n_{j+1}}}\ge \frac 1{2^{j+1}}.
$$
Let $(Z_j)$ be the sequence of subshifts in Lemma \ref{lem:Rothstein} and let 
$Z_\infty$ be their intersection (which has zero topological entropy). Let
$\phi(x)=-d(x,Z_\infty)$. Clearly any maximizing measure for $\phi$ is supported on $Z_\infty$,
and so has integral 0. It follows from the folklore theorem on pressure asymptotics
that $p_\phi(t)\to 0$ as $t\to\infty$. 
We show that $P(t\phi)\ge f(t)$ for all $t\in [t_j,t_{j+1}]$ for all $j\ge j_0$, from
which it follows that $P(t\phi)\ge f(t)$ for all $t\ge t_{j_0}$.

To see this,  we use the variational principle to bound the pressure of $t\phi$ 
from below by $h_T(\mu_j)+\int t\phi\,d\mu_j$ where $\mu_j$ is the (unique) 
measure of maximal entropy of $Z_j$. If $j\ge j_0$ and $t\in [t_j,t_{j+1}]$, then since $\phi\le 0$,
\begin{align*}
P_{\rm top}(t\phi)&\ge P_{\rm top}(t_{j+1}\phi)\ge h_T(\mu_j)+t_{j+1}\int \phi\,d\mu_j\\
&\ge \frac{\log 2}{2^j}-\frac {t_{j+1}}{2^{n_{j+1}}}\\
&\ge \frac{1}{2^{j+1}}=f(t_j)\ge f(t).
\end{align*}
The lower bound for $\int\phi\,d\mu_j$ comes from the second conclusion of Lemma \ref{lem:Rothstein}. 

The main idea in the general case is the same as the above. We briefly describe the 
construction before giving the formal 
proof. The additional complications that we deal with
are (1) we wish to consider pressure functions with non-zero asymptotes; and (2) the ambient shift
$X$ may be a shift of finite type rather than a full shift. 

We find a pair of words, $u$ and $v$ in $\cL(X)$, with the same length and the same first and last letters.
We then use Lemma \ref{lem:Yconstruction} to build a subshift $\tilde X$ of $X$ of entropy exactly $b$ with the property that
the word $u$ appears periodically in points of $\tilde X$ (with some period $l$). We then build a new subshift $X_\infty$ of $X$
by taking the points of $\tilde X$ and replacing the periodic sequence of $u$'s with a Rothstein-coded 
sequence of $u$'s and $v$'s. The function $\phi$ is then $-d(x,X_\infty)$ as in the example above. 

The set $X_\infty$ is the intersection of 
the infinite sequence of $X_j$'s, where the topological entropy can be simply computed. 
Taking $(\mu_j)$ to be the sequence of measures of maximal entropy on $X_j$ (these are unique
although it is not needed for the proof), we have
$p_\phi(t)\ge h(\mu_j)+t\int\phi\,d\mu_j$. Again, provided the sequence $(n_j)$ is suitably constructed,
we can ensure that $p_\phi(t)\ge f(t)$ for all sufficiently large $t$.
\color{black}

\begin{proof}
We observe that for any function $\phi$, $p_{a+\phi}(t)=at+p_\phi(t)$. Hence it suffices to
show that if $f$ is any convex function asymptotic to a constant $b$ such that
$0\le b<h_\text{top}(X)$ then there exists a H\"older potential $\phi$ such that $p_\phi$
is asymptotic to $b$ and $p_\phi(t)>f(t)$ for all sufficiently large $t$.

Since $h_\text{top}(X)>0$, we may pick two distinct words $u,v$ in $\cL(X)$ of the same length
with the same first and last symbols. Let $\ell$ and $Y$ be as constructed in Lemma
\ref{lem:Yconstruction} based on the word $u$,
so that $Y$ is a subshift of topological entropy $b\ell$ whose alphabet
consists of a subset of $\cL_\ell(X)$ with each word in the alphabet starting with $u$.

We next construct a suitable Rothstein shift. Since we assume that $f(t)$ is asymptotic
to the constant $b$, we may pick an increasing
sequence of real numbers $(t_j)$ such that for $t\ge t_j$ we have $f(t)<b+\ell^{-1}2^{-(j+2)}$.
Set $n_j=\lceil \log_2(\ell t_j)\rceil+j$. 
By Lemma \ref{lem:Rothstein} there is a nested family of subshifts
$(Z_j)$ on the alphabet $\{0,1\}$ such that $h_\text{top}(Z_j)=2^{-j}\log 2$ so that $Z_\infty:=\bigcap Z_j$
has $h_\text{top}(Z_\infty)=0$.
Further $d(z,Z_\infty)\le 2^{-n_{j+1}}$ for all $z\in Z_j$. Denote by $\lambda_j$
the measure of maximal entropy for $Z_j$. 

We now let $\bar Y_j=Y\times Z_j\times \{0,\ldots,\ell-1\}$
and consider the constant height $\ell$ suspension over the product $Y\times Z_j$.
That is we consider the map
\begin{equation}\label{eq:susp}
\bar T(y,z,i)=\begin{cases}(Ty,Tz,0)&\text{if $i=\ell-1$;}\\
(y,z,i+1)&\text{otherwise}.
\end{cases}
\end{equation}
Clearly $T\times T$ acting on $Y\times Z_j$ has topological entropy $b\ell+2^{-j}\log 2$
so that $h_\text{top}(\bar Y_j)=b+\log 2/(2^j\ell)$ and $h_\text{top}(\bar Y_\infty)=b$.
Let $\bar\nu_j$ be the invariant measure $\zeta\times \lambda_j\times c$ where $\zeta$ is the measure
of maximal entropy on $Y$ and $c$ is the normalized counting measure on $\{0,\ldots,\ell-1\}$
so that $\bar\nu_j$ is supported on $\bar Y_j$.
Notice that $\bar T^\ell(y,z,i)=(Ty,Tz,i)$ so that $h_{\bar T^l}(\bar\nu_j)=h_T(\zeta)+h_T(\nu_j)
=b\ell+2^{-j}\log 2$. It follows that $h_{\bar T}(\bar\nu_j)=b+\log 2/( 2^j\ell)$ so that
$\bar\nu_j$ is a measure of maximal entropy on $\bar Y_j$.

If the alphabet of $Y$ is $\cW=\{uw_0,\ldots,uw_{N-1}\}$, it is convenient to enlarge it
to $\overline \cW=\cW\cup\{vw_0,\ldots,vw_{N-1}\}$.
We now define a map $\Phi$ from the $\bar Y_j$'s to $X$ by requiring that
$\Phi(y,z,i)_n=\Phi(\bar T^n(y,z,i))_0$ and specifying
$\Phi(y,z,i)_0$:
$$
\Phi(y,z,i)_0=\begin{cases}(y_0)_i&\text{if $i\ge |u|$ or $z_0=0$};\\
v_i&\text{if $i<|u|$ and $z_0=1$.}
\end{cases}
$$
A simple description of $\Phi$ is that it concatenates the words forming
the symbols in $y$, replacing the initial $u$'s with $v$'s in those coordinates
where $z$ has a 1. To see that $\Phi(y,z,i)$ lies in $X$, recall that an arbitrary
finite concatenation of elements of $\cW$ lies in $\cL(X)$.
Since $v\in\cL(X)$ and has the same first and last symbols as $u$, replacing
any number of $u$'s by $v$'s in an element of $\cL(X)$ gives another element
of $\cL(X)$. It follows that arbitrary finite concatenations of elements of $\overline\cW$
lie in $\cL(X)$.

We let $X_j=\Phi(\bar Y_j)$ for all $1\le j\le\infty$
and $\mu_j=\bar\nu_j\circ\Phi^{-1}$. We now claim that $\Phi$ is
at most $\ell$-to-one and hence entropy preserving.
We will write $\Phi$ as a composition $\kappa\circ\Psi$.
Firstly, we introduce the space $\overline\cW^\Z\times \{0,\ldots,\ell-1\}$ with the suspension map
$\bar T$ defined exactly like \eqref{eq:susp}.
Then $\Psi:\bar Y_j\to \overline\cW^\Z\times \{0,\ldots,\ell-1\}$
is defined by $\Psi(y,z,i)=(s(y_n,z_n)_{n\in\Z},i)$, where
$$
s(y,z)_n=\begin{cases}
y_n&\text{if $z_n=0$};\\
v*(y_n)_{|u|}^{\ell-1}
&\text{if $z_n=1$,}
\end{cases}
$$
where $*$ denotes concatenation and $y_i^j$ denotes the word $y_iy_{i+1}\ldots y_j$. 
That is, $s$ replaces the initial $u$'s with $v$'s in those coordinates where $z_n=1$.
This mapping is invertible as one can infer the sequence $z$ from looking at $s(y,z)$.
Secondly, the map $\kappa\colon \overline\cW^\Z\times \{0,\ldots,\ell-1\}\to X$ is a concatenation map.
$\kappa(y,i)_n=\kappa(\bar T^n(y,i))_0$ and
$\kappa(y,i)_0=(y_0)_i$. The map $\kappa$ is at most $\ell$-to-1 as given
$\kappa(y,i)$ and the value of $i$, one may recover $y$.
Since $\Phi$ is entropy-preserving, we see
\begin{equation}\label{eq:YxZentropy}
h_T(\mu_j)=h_{\bar T}(\nu_j)=b+\frac{\log 2}{2^j\ell}\text{ for all $1\le j\le\infty$.}
\end{equation}

To finish the argument, we define $\phi(x)=-\dist(x,X_\infty)$. Clearly any maximizing measure $\mu_\infty$
is supported on $X_\infty$ so that $\int\phi\,d\mu_\infty=0$ and $h_T(\mu_\infty)=
h_\text{top}(X_\infty)=h_\text{top}(\bar Y_\infty)=b$. By the folklore theorem on pressure asymptotics, this ensures
that $p_\phi(t)$ is asymptotic to $b$ as required.

Let $x\in X_j$ so that $x=\Phi(y,z,i)$ with $y\in Y$, $z\in Z_j$ and $0\le i<\ell$.
By Lemma \ref{lem:Rothstein}, there exists $z'\in Z_\infty$ with $d(z,z')<2^{-n_{j+1}}$.
By definition $\Phi(y,z',i)\in X_\infty$. From the definition of $\Phi$, one sees that
$x_n=\Phi(y,z,i)_n=\Phi(y,z',i)_n$ for all $|n|<n_{j+1}\ell$ so that
\begin{equation}\label{eq:phi_on_Zj}
d(x,X_\infty)\le 2^{-n_{j+1}\ell}\text{ for all $x\in X_j$}.
\end{equation}

We claim that for $t_j\le t\le t_{j+1}$, $P(t\phi)\ge f(t)$.
As before, we rely on the Variational Principle and just show that $h(\mu_j)+t\int \phi\,d\mu_j\ge f(t)$. 
Combining \eqref{eq:YxZentropy} and \eqref{eq:phi_on_Zj} for $t_j\le t\le t_{j+1}$ we obtain
  \begin{align*}
    P(t\phi) & \ge h_T(\mu_j)+t\int\phi\,d\mu_j \\
    &= b+\frac{1}{\ell 2^j}-t\int \dist(x,X_\infty)\, d\mu_j\\
    &\ge b+\frac{\log 2}{\ell 2^j}-\frac{t_{j+1}}{\ell t_{j+1}2^{j+2}}\\
    &\ge b+\frac{1}{\ell 2^{j+2}}\\
    &\ge f(t).
  \end{align*}
\end{proof}

\section{Genericity of upper bounds}\label{sec:genericity}
First we give examples of two potentials on the full shift on two symbols for which the behaviour of the  pressure functions is very
different even though they share the same maximizing measure.
\begin{example}
  There is a potential $\phi$ on the full shift $\{0,1\}^\Z$ satisfying
  \begin{itemize}
    \item the unique ground state of $\phi$ is a point-mass measure at $\bar{0}$;
    \item $p_{\phi}(t)-\ell_\infty(t) \sim e^{-t}$.
  \end{itemize}
\end{example}

\begin{proof}
Let $X=\{0,1\}^\Z$ and let $\phi(x)=-x_0$.
Then the unique maximizing measure is $\delta_0$, the measure supported on the fixed point
$\bar 0$ so that $h(\delta_0)=0$ and $\int \phi\,d\delta_0=0$. Hence $\ell_\infty(t)=0$. We see from Example \ref{ex:optimal_bnd} that
$P_{\rm top}(t\phi)=\log(1+e^{-t})=e^{-t}+O(e^{-2t})$. In particular, $p_\phi(t)$
approaches $\ell_\infty(t)$ exponentially fast.
\end{proof}{}

\begin{example}\label{ex:no_upper_bnd}
  There is a potential $\phi$ on the full shift $\{0,1\}^\Z$ satisfying
  \begin{itemize}
    \item the unique ground state of $\phi$ is a point-mass measure at $\bar{0}$;
    \item $p_{\phi}(t)-\ell_\infty(t)\gtrsim \log\log t/\log t$.
  \end{itemize}
\end{example}
\begin{proof}
Let $X=\{0,1\}^\Z$ and let $$S=\{x\in X\colon x \text{ has at most one 1 symbol}\}.$$ Then $S$ is a countable closed
subshift of $X$, sometimes known as the ``sunny side up" system. Let $\phi(x)=-d(x,S)$. Clearly the only
invariant measure with support lying in $S$ is $\delta_0$, so that $\ell_\infty(t)=0$ again. We show,
however, that in this case $p_\phi(t)$ does not approach $\ell_\infty$ exponentially fast. Indeed
let $\mu_n$ be the measure on $X$ where gaps are uniform and equally likely in the range $2n,\ldots,3n$.
Then a calculation using Abramov's theorem shows that $h(\mu_n)=(\log(n+1))/(\frac {5n}2)$ (one chooses
between $n+1$ equally likely choices on average once every $\frac{5n}2$ steps).
For $\mu_n$-a.e.\ $x$, $x_{-(n-1)},\ldots x_{n-1}$ contains at most one 1, so that $d(x,S)\le 2^{-n}$.
It follows that $\int \phi\,d\mu_n\ge -2^{-n}$. We now have the lower bound for $p_\phi(t)=P(t\phi)\ge
\max_n (\frac 25\log n/n-t2^{-n})$. Taking $n=\lceil \log_2 t+\log_2(\log_2 t)\rceil$, we see that for large $t$,
$p_\phi(t)\ge \frac 15\log(\log_2 t)/\log_2 t-1/(\log_2 t)$. So that $p_\phi(t)-\ell_\infty(t)$ converges
to zero \emph{much} slower than exponentially.
\end{proof}

Before establishing Theorem \ref{thm:gen_upper_bnd}, we prove the following lemma which is a version
of a lemma of Yuan and Hunt \cite[Proposition 4.1]{YH}. Let $X$ be a mixing shift of finite
type. For $\alpha<1$, let $d_\alpha(x,y)=\alpha^{\min\{|n|\colon x_n\ne y_n\}}$ and
$\|\phi\|_\alpha=\max(\|\phi\|_\infty,
\sup_{x\ne y}|\phi(x)-\phi(y)|/d_\alpha(x,y))$.
We then define the H\"older functions to be $\cH_\alpha=\{\phi\colon \|\phi\|_\alpha<\infty\}$.

\begin{lem}\label{lem:YH}
Let $X$ be a shift of finite type and let $p$ be a periodic point with periodic orbit $\cO(p)$. 
Let $\mu_p$ denote the invariant measure supported on $\cO(p)$. 
Let $d_p(x)=d_\alpha(x,\cO(p))$.
There exists $r>0$ such that 
if $\|\psi\|_\alpha<r$, then the function $\phi=-d_p+\psi$ satisfies
$\int\phi\,d\mu_p<\int\phi\,d\mu$ for all $T$-invariant
measures $\mu\ne \mu_p$. That is, $\mu_p$ is the unique minimizing measure for all $\phi$ 
of the form $-d_p+\psi$ with $\|\psi\|_\alpha<r$.
\end{lem}

\begin{proof}
Let $k$ be the period of $p$. Let $\gamma=\min\{d_\alpha(T^ip,T^jp)\colon 0\le i<j<p\}$ 
and $\eta=\alpha^{k}\gamma/2$.
This has the property that if $d_\alpha(x,T^ip)\le \eta$, then for
$n=0,\ldots,k-1$, the closest point of $\cO(p)$ to $T^nx$
is $T^{n+i}p$ and 
$d_\alpha(T^nx,T^{n+i}p)\le \alpha^{-n}d(x,T^ip)$. In this case, 
we say $x$ \emph{follows} $\cO(p)$ for a cycle. If $d_\alpha(T^{jk}x,\cO(p))<\eta$ for $j=0,\ldots m-1$,
but $d_\alpha(T^{mk}x,\cO(p))\ge\eta$, we say $x$ follows $\cO(p)$ for $m$ cycles. If $d_\alpha(T^{mk}x,\cO(p))<\eta$
for all $m\ge 0$, then $T^nx$ approaches $\cO(p)$ exponentially by the shadowing property of hyperbolic dynamical systems. 

Let $r=\frac 14(1-\alpha)\alpha^{k+1}$
and let $\phi=-d_p+\psi$ with $\|\psi\|_\alpha<r$. Let $\ell=\int\phi\,d\mu_p$.
Then we have $\ell=\int\psi\,d\mu_p=\frac 1k\sum_{x\in\cO(p)}\psi(x)$, so that $\ell\ge -r$. 
Let $x\in X$ be such that $d_\alpha(T^nx,\cO(p))\ge \eta$ for infinitely many $n$.
We consider two cases: either $d(x,\cO(p))\ge\eta$ or $x$ follows $\cO(p)$ for $m$ cycles for some $m\ge 1$.
In the first case, we estimate 
$$
\phi(x)-\ell\le d_p(x)+2r=-\eta+2r<-\tfrac \eta 2.
$$
and in the second case, we estimate $S_{km+1}\phi(x)-(km+1)\ell$. 
Let $T^ip$ be the closest point of $\cO(p)$ to $x$. We have
\begin{align*}
S_{km+1}(\phi-\ell)(x)&=
S_{km}\phi(x)-S_{km}\phi(T^ip)+\phi(T^{km}x)-\ell\\
&\le \sum_{j=0}^{km-1}|\phi(T^jx)-\phi(T^j(T^ip))|-\eta+2r.
\end{align*}
The last $k$ terms of the sum may be estimated above using the Lipschitz
condition by $r(\eta+\eta/\alpha+\ldots + \eta/\alpha^{k-1})$.
The previous terms may be estimated above using the hyperbolic shadowing lemma and the Lipschitz condition by 
$2r(\eta+\eta\alpha+\eta\alpha^2+\ldots)$. 
Combining the estimates, we see 
\begin{align*}
S_{km+1}(\phi-\ell)(x)&\le 2r\alpha^{-(k-1)}/(1-\alpha)-\eta+2r\\
&<2r\frac{\alpha^{-(k+1)}}{1-\alpha}-\eta=-\tfrac\eta 2.
\end{align*}
Thus if $\mu$ is any ergodic invariant measure not supported on $\cO(p)$, for $\mu$-a.e. $x$,
$x$ is not asymptotic to $\cO(p)$ so the orbit may be split into an infinite number of segments of
the two types described above, with an orbit sum of of $\phi-\ell$ being at most $-\eta/2$ on each segment. 
It follows, using Birkhoff's theorem, that $\int\phi\,d\mu-\int\phi\,d\mu_p=\int (\phi-\ell)\,d\mu<0$,
so that $\mu_p$ is the unique minimizing measure for $\phi$ as claimed. 
\end{proof}

\color{black}
We conclude the paper with the proof of Theorem \ref{thm:gen_upper_bnd}, which we restate for convenience.
\begin{thm*}
Let $X$ be a mixing shift of finite type and let $\cH_\alpha$
be as above.
Then there is a dense open subset $\cU$ of potentials in $\cH_\alpha$ such that for all $\phi\in \cU$,
there exist $C_1>0$ and $C_2>0$ such that $p_\phi(t)-\ell_\infty(t)\le C_1e^{-C_2t}$ for all $t$.
\end{thm*}

\begin{proof}
Note that $\cH_\alpha$ is the class of functions that is Lipschitz with respect to $d_\alpha$ where
$d_\alpha(x,y)=\alpha^{\min\{|n|\colon x_n\ne y_n\}}$. Let $\mathsf{Per}$ denote
the collection of periodic points in $X$.
For a point $p\in\mathsf{Per}$, let $\cH_\alpha(p)$ denote the set potentials $\phi$ in $\cH_\alpha$ such that
$\int \phi\,d\mu_p\le \int \phi\,d\nu$ for
all invariant measures $\nu\ne \mu_p$ (where $\mu_p$ denotes the unique invariant measure
supported on the orbit of $p$).

Recall that by a theorem of Contreras \cite{Co}, $\bigcup_{p\in\mathsf{Per}}\cH_\alpha(p)$ contains a dense open
subset of $\cH_\alpha$. For a periodic point $p$, we let 
$d_p(x)=\min_jd_\alpha(x,\cO(p))$ as above.

The set $\cU$ that we consider is defined by
$$
\cU=\bigcup_{p\in\mathsf{Per}}\{\phi-ad_p\colon \phi\in \cH_\alpha(p), a>0\}.
$$
Since $\bigcup_{p\in\mathsf{Per}} \cH_\alpha(p)$ is in the closure of $\cU$, we see that $\cU$ is dense.
To show that $\cU$ is open, let $p$ be a periodic point, 
$\phi\in \cH_\alpha(p)$ and $a>0$. We claim that if $\|\psi\|_\alpha<\frac a2r$, where $r$
is as in the statement of Lemma \ref{lem:YH}, then $\phi-ad_p+\psi\in \cU$.
To see this, we note that $\phi\in\cH_\alpha(p)$ and $-\frac a2d_p+\psi\in\cH_\alpha(p)$ by Lemma \ref{lem:YH},
so $\phi-\frac a2d_p+\psi\in \cH_\alpha(p)$. Hence $\phi-ad_p+\psi=(\phi-\frac a2d_p+\psi)-\frac a2d_p\in \cU$.


We now let $\phi=\psi-ad_p$ where $p\in\mathsf{Per}$, $\psi\in \cH_\alpha(p)$ and $a>0$. Since $a>0$, we
see that $\int\phi\,d\nu<\int\phi\,d\mu_p$ for all invariant measures $\nu\ne\mu_p$. It follows that
$\ell_\infty(t)$, the tangent line to $p_\phi(t)$ at infinity, is given by $\ell_\infty(t)=\beta t$
where $\beta=\int\phi\,d\mu_p=\int\psi\,d\mu_p$.
By the ``Ma\~n\'e lemma" (see \cite[Theorem 4.7]{Jenkinson} for a version that applies in our setting),
there exists a H\"older continuous function
$\tilde\psi$ cohomologous to $\psi$ such that
$\tilde\psi(x)\le\beta$ for all $x\in X$. Let $\tilde\phi=\tilde\psi-ad_p$
and notice that $t\tilde\phi$ is cohomologous to $t\phi$ for all $t$ so that
$$p_\phi(t)=P_{\rm top}(t\tilde\phi)=P_{\rm top}(t\tilde\psi-ad_p)\le P_{\rm top}(\beta t-atd_p)=\beta t+P_{\rm top}(-atd_p).$$
Hence it suffices to show that $P_{\rm top}(-atd_p)$ decreases exponentially to 0.
Since $\mu_p$ is the unique maximizing measure for $-atd_p$,
the folklore theorem on pressure asymptotics shows that $P_{\rm top}(-atd_p)\ge 0$ for all $t$ and it converges
to 0 as $t\to\infty$.
We give a crude bound showing the exponential convergence.
First, we enlarge the space $X$ to the full shift $\bar X$ whose
alphabet is the alphabet $\cA$ of $X$. Then define a potential $\bar f$ on $\bar X$
by
$$
\bar f(x)=\begin{cases}-ad_p(x)&\text{if $x_{-1}x_0$ is legal in $X$;}\\
-\infty&\text{otherwise;}
\end{cases}
$$
so that $P_{\rm top}(X,-atd_p)=P_{\rm top}(\bar X,t\bar f)$. Let $k$ be the period of $p$. We then observe that
$\bar f(x)\le \bar g(x)$ for all $x$ where
$$
\bar g(x)=\begin{cases}
-a/2^{k}&\text{if $x_0\ne x_{-k}$;}\\
0&\text{otherwise.}
\end{cases}
$$
Hence it follows that
$P_{\rm top}(\bar X,t\bar f)\le P_{\rm top}({\bar X},t\bar g)$ for all $t>0$.
Since the maximizing measures for $\bar g$ are precisely those measures supported on period $k$ orbits,
we see by the folklore theorem on pressure asymptotics that $P_{\rm top}(\bar X,t\bar g)\to 0$ as $t\to\infty$.
Summing $e^{S_{kn}t\bar g(x)}$
over words of length $nk$ shows that $P_{\rm top}({\bar X},t\bar g)=\log(1+(|\cA|-1)e^{-at/2^k}))$.
As noted above
$P_{\rm top}(X,t\phi)=\beta t+P_{\rm top}(X,-atd_p)\le \beta t+P_{\rm top}({\bar X},t\bar g)$. Combining this with the
inequality that we just derived shows that
$p_\phi(t)$ approaches the
asymptote $\ell_\infty$ exponentially fast.
 \end{proof}

\end{document}